\documentclass[10pt]{amsart}
\usepackage{enumitem}
\usepackage{graphicx}
\usepackage{marvosym}
\usepackage{amssymb}

\usepackage{upref}
\usepackage{hyperref}

\newtheorem{definition}{Definition}[section]
\newtheorem{theorem}{Theorem}[section]
\newtheorem{lemma}{Lemma}[section]

\newtheorem*{maintheorem*}{Main Theorem}
\allowdisplaybreaks
\numberwithin{equation}{section}

\newcommand{\et}{\wedge}
\newcommand{\vel}{\vee}
\newcommand{\weakstar}{\overset{\star}\rightharpoonup}

\newcommand{\norm}[1]{\left\Vert#1\right\Vert}
\newcommand{\abs}[1]{\left|#1\right|}
\newcommand{\Set}[1]{\left\{#1\right\}}

\newcommand{\eps}{\varepsilon}

\newcommand{\pt}{\partial_t}
\newcommand{\px}{\partial_x }
\newcommand{\pxx}{\partial_{xx}^2}
\newcommand{\pxxx}{\partial_{xxx}^3}

\newcommand{\ptxx}{\partial_{txx}^3}

\newcommand{\wlim}[1]{\overline{#1}}

\renewcommand{\i}{\ifmmode\mathit{\mathchar"7010 }\else\char"10 \fi}
\renewcommand{\j}{\ifmmode\mathit{\mathchar"7011 }\else\char"11 \fi}
\newcommand{\R}{\mathbb{R}}

\newcommand{\Do}{(0,T)\times \R}
\newcommand{\hDo}{[0,T)\times \R}
\newcommand{\cDo}{[0,T]\times \R}

\newcommand{\loc}{\mathrm{loc}}

\newcommand{\supp}{\mathrm{supp}\,}

\newcommand{\dx}{\, dx}
\newcommand{\dy}{\, dy}
\newcommand{\dt}{\, dt}

\newcommand{\seq}[1]{\left\{#1\right\}}

\newcommand{\test}{\varphi}

\newcommand{\Z}{\mathbb{Z}}

\newcommand{\weak}{\rightharpoonup}
\newcommand{\Dx}{{\Delta x}}
\newcommand{\Dt}{{\Delta t}}

\newcommand{\dto}{\downarrow}

\newcommand{\Dm}{D_{-}}
\newcommand{\Dp}{D_{+}}

\newcommand{\Dpm}{D_\pm}
\newcommand{\Dmp}{D_\mp}

\newcommand{\tqnj}{\tilde{q}_{j}^n}
\newcommand{\tqnjpe}{\tilde{q}_{j+1}^n}

\newcommand{\tunjph}{\tilde{u}_{j+1/2}^n}
\newcommand{\tunjmh}{\tilde{u}_{j-1/2}^n}

\newcommand{\tPnj}{\tilde{P}_{j}^n}

\newcommand{\unjph}{u_{j+1/2}^n}
\newcommand{\unjmh}{u_{j-1/2}^n}

\newcommand{\Pnj}{P_{j}^n}

\newcommand{\Pj}{P_{j}}
\newcommand{\Pdx}{P_{\Dx}}

\newcommand{\qj}{q_{j}}

\newcommand{\qnj}{q_{j}^n}
\newcommand{\qnjp}{q_{j+1}^n}
\newcommand{\qnjm}{q_{j-1}^n}

\newcommand{\ujph}{u_{j+1/2}}
\newcommand{\ujmh}{u_{j-1/2}}

\newcommand{\DT}{D_{+}^t}
\newcommand{\weakto}{\rightharpoonup}
\newcommand{\el}[1]{\ell^{#1}} 
\newcommand{\fm}{f^M}
\newcommand{\udx}{u_{\Dx}}
\newcommand{\qdx}{q_{\Dx}}
\newcommand{\ee}{\eta_{\eps}}

\newcounter{asnr}

{\ifnum\value{asnr}=0 \stepcounter{asnr} 
  \begin{enumerate}[label=\textbf{A}.\arabic{enumi}]
    \else
    \begin{enumerate}[label=\textbf{A}.\arabic{enumi},resume] \fi}
{\end{enumerate}}

\newcounter{defnr}

\newenvironment{Definitions} %
{\ifnum\value{defnr}=0 \stepcounter{defnr} 
  \begin{enumerate}[label=(\textbf{D}.\textbf{\arabic{enumi}})]
    \else
    \begin{enumerate}[label=(\textbf{D}.\textbf{\arabic{enumi}}),resume] \fi}
{\end{enumerate}}

\begin{document}

\title[An explicit scheme for the Camassa-Holm equation]
{An explicit finite difference scheme \\ for the Camassa-Holm equation}

\date{\today}

\author[G. M. Coclite]{G. M. Coclite} 
\address[Giuseppe Maria Coclite]{\newline Dipartimento di Matematica \newline
Universit\`a degli Studi di Bari \newline
Via E. Orabona 4\newline
70125 Bari, Italy}
\email[]{coclitegm@dm.uniba.it}

\author[K. H. Karlsen]{K. H. Karlsen} 
\address[Kenneth H. Karlsen]{\newline Centre of Mathematics for Applications
  (CMA) \newline University of Oslo\newline P.O. Box 1053,
  Blindern\newline N--0316 Oslo, Norway}
\email[]{kennethk@math.uio.no}
\urladdr{http://folk.uio.no/kennethk/}

\author[N. H. Risebro]{N. H. Risebro} 
\address[Nils Henrik Risebro]{\newline Centre of Mathematics for Applications (CMA)\newline
  University of Oslo\newline P.O. Box 1053, Blindern\newline N--0316
  Oslo, Norway} \email[]{nilshr@math.uio.no}
\urladdr{http://folk.uio.no/nilshr/}

\subjclass[2000]{35G25, 35L05, 65M06, 65M12}






\keywords{Camassa-Holm equation, hyperbolic-elliptic system, weak solution, 
finite difference scheme, convergence}  

\thanks{This work was supported by the Research Council of Norway through 
the project WaveMaker and an Outstanding Young Investigators Award of K. H. Karlsen.}   

\begin{abstract}
We put forward and analyze an explicit finite difference scheme for the 
Camassa-Holm shallow water equation that can handle general $H^1$ initial data and thus 
peakon-antipeakon interactions. Assuming a specified condition restricting the 
time step in terms of the spatial discretization parameter, we prove that the 
difference scheme converges strongly in $H^1$ towards a dissipative weak 
solution of Camassa-Holm equation.
\end{abstract}

\maketitle

\tableofcontents

\section{Introduction}\label{seq:intro} 
In this paper we present and analyze an explicit finite difference scheme for
the Camassa-Holm partial differential equation \cite{Camassa:1993zr}
\begin{equation}\label{eq:u}
        \pt u-\ptxx u+ 3u \px u=2 \px u  \pxx u + u \pxxx u, 
        \qquad (t,x)\in \Do,
\end{equation}
which we augment with an initial condition:
\begin{equation}\label{ass:init}
        u|_{t=0}=u_0, \qquad u_0\in H^1(\R), u_0\neq 0.
\end{equation}
Rewriting equation \eqref{eq:u} as
$$
(1-\pxx) \left[\pt u + u \px u\right] + \px \left(u^2+\frac{1}{2}
  (\px u)^2\right)=0,
$$
we see that (for smooth solutions) \eqref{eq:u} is equivalent to
the elliptic-hyperbolic system
\begin{equation}\label{eq:uprime}
        \pt u +u \px u + \px P=0,
        \qquad - \pxx P+P=u^2+\frac{1}{2}(\px u)^2.
\end{equation}
Recalling that $e^{-|x|}/2$ is the Green's function of the operator
$1-\pxx$, \eqref{eq:uprime} can be written as
\begin{equation}\label{eq:conslaw}
        \pt u + \px F(u,\px u)=0, \quad 
        F(u,\px u)= \frac12\left[u^2 + e^{-\abs{x}}\star
        \left(u^2 +\frac{1}{2} (\px u)^2\right)\right],
\end{equation}
which can be viewed as a conservation law with nonlocal flux function.
In this paper the relevant formulation of the Camassa-Holm equation
\eqref{eq:u} is the one provided by the hyperbolic-elliptic system
\eqref{eq:uprime} or \eqref{eq:conslaw}.

The Camassa-Holm equation can be viewed as a model for the propagation
of unidirectional shallow water waves
\cite{Camassa:1993zr,Johnson:2002bs}; it is a member of the class of
weakly nonlinear and weakly dispersive shallow water models, a class
which already contains the Korteweg-de Vries (KdV) and
Benjamin-Bona-Mahony (BBM) equations. In another interpretation the
Camassa-Holm equation models finite length, small-amplitude radial
deformation waves in cylindrical compressible hyperelastic rods
\cite{Dai:1998fu}. It arises also in the context of differential
geometry as an equation for geodesics of the $H^1$-metric on the
diffeomorphism group, see for example
\cite{Constantin:2002ij,Constantin:2003uw,Holm:1998hc,Misioek:1998tg}.
The Camassa-Holm equation possesses several striking properties such
as an inifinite number of conserved integrals, a bi-Hamiltonian
structure, and complete integrability
\cite{Beals:2000fk,Camassa:1993zr,Constantin:1999ly,Constantin:2001ve,Fuchssteiner:1981fk}.
Moreover, it enjoys an infinite number of non-smooth solitary wave
solutions, called peakons, which are weak solutions of
\eqref{eq:conslaw}.
 
From a mathematical point of view the Camassa-Holm equation has by now
become rather well-studied. While it is impossible to give a complete
overview of the mathematical literature, we shall here mention a few
typical results, starting with the local(-in-time) existence results
in \cite{Constantin:1998uq,LiOlver:2000,Rodriguez-Blanco:2001fk} and
those using Besov spaces in \cite{Danchin:2003bh,Danchin:2001dq}.  It
is well-known that global solutions do not exist and wave-breaking
occurs \cite{Camassa:1993zr}.  Wave-breaking means that the solution
itself stays bounded while the spatial derivative becomes unbounded in
finite time.

In view of what we have said so far (peakon solutions/wave-breaking)
it is clear that a theory based on weak solutions is essential.  In
the literature there are a number of results on (dissipative and
conservative) weak solutions of the Camassa-Holm equation, see
\cite{Bressan:2006oq,Bressan:2006nx,Bressan:2006wd,Constantin:1998kx,Coclite:2006ny,Constantin:2000zr,Holden:2006rr,Xin:2000qf,Xin:2002ys}
and the references cited therein.  In this paper we are interested
specifically in the class of \textit{dissipative weak solutions}
studied by Xin and Zhang \cite{Xin:2000qf,Xin:2002ys}. Their results
show, among other things, that there exists a global dissipative weak
solution of \eqref{eq:u}-\eqref{ass:init} for any $H^1$ initial data
$u_0$ (peakon-antipeakon interactions are covered).  These solutions
are global in the sense that they are defined past the blow-up time
(wave-breaking).  More precisely, suppose $u_0 \in H^1(\R)$. Then
there exists a global weak (distributional) solution $u\in L^\infty(0,
T ; H^{1}(\R))$ of \eqref{eq:u} satisfying the following properties:
$t\mapsto \norm{u(t,\cdot)}_{H^1(\R)}$ is non-increasing; $\px u \in
L^p_{\loc}(\R_+\times \R)$, $p<3$; 
\begin{equation}
  \px u (t,x) \le \frac{2}{t} +
  C \norm{u_0}_{H^1(\R)}
  ,\quad \text{for} \quad
  t>0, \label{eq:ent1}
\end{equation}
for some positive constant $C$. This last item presumably singles out
a unique weak solution. As an example of how this may work we consider
the ``peakon-antipeakon'' solution given by
\begin{equation}
u(t,x)=\tanh(t-1)\left(e^{-\abs{x-y(t-1)}}-e^{-\abs{x+y(t-1)}}\right),
\quad
y(t)=\log\left(\cosh(t)\right).\label{eq:pap}
\end{equation}
This formula represents a peakon ($e^{-\abs{x+y}}$) colliding with an
antipeakon ($-e^{-\abs{x-y}}$) at $x=0$ and $t=1$. Note that
$u(1,x)=0$. How this solution
is extended to $t>1$ depends on which solution concept we adopt. If we
use the formula \eqref{eq:pap} also for $t>1$ we get the conservative
solution for which $\norm{u(t,\cdot)}_{H^1(\R)}$ is constant for
almost all $t$. We can
also extend the solution by defining $u(t,x)=0$ for $t>1$. Obviously,
the ``entropy condition'' \eqref{eq:ent1} will only be satisfied for
this dissipative solution. 

Let us now turn to the topic of the present paper, which is the design
and analysis of numerical schemes.  The first numerical results for
the Camassa-Holm equation are presented in \cite{Camassa:1994pd} using
a pseudo-spectral scheme. Numerical simulations with pesudo-spectral
schemes are also reported in \cite{Fringer:2001kl,Holm:2003qe}.
Numerical schemes based on multipeakons (thereby exploiting the
Hamiltonian structure of the Camassa-Holm equation) are examined in
\cite{Camassa:2003lq,Camassa:2005rr,Camassa:2006nx}.  In
\cite{Holden:2006oq}, the authors prove that the multipeakon algorithm
from \cite{Camassa:2005rr,Camassa:2006nx} converges to the solution of
the Camassa-Holm equation as the number of peakons tends to infinity.
This convergence result applies to the specific situation where the
initial function $u_0\in H^1$ is such that $(1 -\pxx) u_0$ is a
positive measure. For the same class of initial data, in
\cite{Holden:2005pr} the authors prove that a semi-discrete finite
difference scheme based on the variable $m:=(1 -\pxx) u$ converges
strongly in $H^1$ to the weak solution identified in
\cite{Constantin:1998kx,Constantin:2000zr}. In \cite{Kalisch:2006kx},
the authors establish error estimates for a spectral projection scheme
for smooth solutions. In a different direction, an adaptive
high-resolution finite volume scheme is developed and used in
\cite{Artebrant:2006qy}.  The local discontinuous Galerkin method is
adapted to the Camassa-Holm equation in \cite{Xu:2007lr}. Although
this work does not provide a rigorous convergence result for general
(non-smooth) solutions, they show that the discrete total energy is
nonincreasing in time, thereby suggesting that the approximate
solutions are of dissipative nature. Besides, they establish an error
estimate for smooth solutions.  Finally, multi-symplectic schemes
possessing good conservative properties are suggested and demonstrated
in the recent work \cite{Cochen:2007fk}.

It seems rather difficult to construct numerical schemes for which one
can prove the convergence to a (non-smooth) solution of the Camassa-Holm
equation. This statement is particularly accurate in the case of 
general $H^1$ initial data and peakon-antipeakon interactions. 
Indeed, in this context we are only aware of the recent work \cite{Coclite:2006kx} in which 
we prove convergence of a tailored semi-discrete difference scheme 
to a dissipative weak solution.  Before we can outline this scheme, let us discretize the 
spatial domain $\R$ by specifying the mesh points 
$x_j=j\Dx$, $x_{j+1/2}=(j+1/2)\Dx$, $j=0,\pm 1, \pm 2,\dots$, where 
$\Dx>0$ is the length between two consecutive mesh points (the spatial discretization 
parameter). Let $D_-$, $D$, and $D_+$ denote the corresponding 
backward, central, and forward difference operators, respectively. 
The scheme proposed in \cite{Coclite:2006kx}, which 
is based on the formulation \eqref{eq:uprime}, reads
\begin{equation}\label{eq:u-disc-semi}
        \begin{split}
                &\frac{d}{dt}\ujph + \left(\ujph\vee 0 \right)\Dm\ujph 
                + \left(\ujph\wedge 0 \right)\Dp\ujph +\Dp\Pj=0,
                \\ & -D_-D_+\Pj+\Pj =\left(\ujph\vee 0 \right)^2
                +\left(\ujmh\wedge 0 \right)^2 
                + \frac12 \left(\Dm\ujph\right)^2,
        \end{split}
\end{equation}
where 
$$
\ujph(t)\approx u(t,x_{j+1/2}), \quad P_j(t)\approx P(t,x_j), \quad 
\text{for $t\ge 0$ and $j\in \Z$}.
$$ 

If we interpret the Camassa-Holm equation \eqref{eq:conslaw} as a ``perturbation" of the inviscid 
Burgers equation, then the $u$-part of \eqref{eq:u-disc-semi} might not come across as a reasonable 
(upwind) difference scheme. On the other hand, as pointed out in \cite{Coclite:2006kx}, the key point is that 
with \eqref{eq:u-disc-semi} the quantity $\qj:=\Dm\ujph$ satisfies a difference 
scheme which contains proper upwinding of the transport term in the equation 
for $q:=\px u$, which reads $\pt q + u \px q + \frac{q^2}{2}+P-u^2=0$. 
Consequently, as is proved in \cite{Coclite:2006kx}, the scheme \eqref{eq:u-disc-semi} 
satisfies a total energy inequality, in which only the $q$-part of the total energy is dissipated (not 
the $u$-part, which is after all continuous). This is the essential starting point for the 
entire convergence analysis in \cite{Coclite:2006kx}.

The ``semi-discrete" equation in \eqref{eq:u-disc-semi} constitutes an infinite system of 
ordinary differential equations which must be solved by some numerical method. 
The main purpose of the present paper is to show
that a fully discrete version of the scheme used in
\cite{Coclite:2006kx} produces a convergent sequence of approximate
solutions, and that the limit is a dissipative weak solution to \eqref{eq:u}. 
The fully discrete version that is analyzed in this paper is based on
replacing the time derivative in \eqref{eq:u-disc-semi} by a forward
difference, i.e.,
$$
\ujph'(t)\to \DT\unjph := \frac{u^{n+1}_{j+1/2}-\unjph}{\Dt},
$$
and evaluating the rest of \eqref{eq:u-disc-semi} at $t^n:=n\Dt$. Now
$\unjph$ should approximate the exact solution $u$ at the 
point $(t^n,x_{j+1/2})$. This gives the fully discrete scheme
\begin{equation}\label{eq:u-disc-fully}
  \begin{split}
    \DT\unjph + \left(\unjph\vee 0\right)\Dm\unjph 
    + \left(\unjph\wedge 0\right)\Dp\unjph + \Dp\Pnj =0\\
    -\Dm\Dp\Pnj + \Pnj=\left(\unjmh\vee 0\right)^2 +
    \left(\unjph\wedge 0\right)^2  
    + \frac{1}{2}\left(\Dm\unjph\right)^2,
  \end{split}
\end{equation}
where $\Pnj$ approximates $P(t^n,x_j)$. As in
\cite{Coclite:2006kx} this is a difference scheme which is tailored so
that it gives an upwind scheme for the equation satisfied by $q:=\px u$.  

The main aim of this paper is prove that the fully discrete (explicit) scheme \eqref{eq:u-disc-fully} 
converges to a dissipative weak solution of the Camassa-Holm equation. 
The starting point of the analysis is a total energy
estimate, showing that the $H^1$ norm of the approximate 
solutions is (almost) nonincreasing in time. To this end, we must 
assume that
\begin{equation}\label{eq:cfl1}
        \Dt=\mathcal{O}(\Dx^2\log(1+\Dx^\theta))
\end{equation}
for some $\theta>0$
as $\Dx\to 0$. This is a very severe condition, and it may seem that
when using this method in practice one should use very small 
time steps. However, this is not a Neumann type stability criterion,
and we do not have blow up if it is violated. Indeed, practical
experiments indicate stability and convergence if $\Dt=\mathcal{O}(\Dx)$. 

By appropriately extending the difference solution
\eqref{eq:u-disc-fully} to a function $\udx(t,x)$ defined at all
points $(t,x)$ in the domain, we prove under condition \eqref{eq:cfl1}
that $\seq{\udx}_{\Dx>0}$ converges strongly in $H^1$ to a dissipative
weak solution of the Camassa-Holm equation
\eqref{eq:u}-\eqref{ass:init}.  Regarding the proof, we adapt the
``renormalization" approach used in \cite{Coclite:2006kx} for the
semi-discrete scheme, but there are several essential deviations and
many parts of the convergence proof are substantially more involved
and/or different.  These differences are mainly due to the fact that
the semi-discrete scheme, when viewed as a fully discrete scheme with
``infinitely small time steps'', has a large and stabilizing numerical
viscosity.  Regarding the fully discrete (explict) scheme
\eqref{eq:u-disc-fully}, to account for this lack of numerical
viscosity the convergence analysis relies heavily on the CFL condition
\eqref{eq:cfl1} and differennt With reference to the differences
between the semi-discrete and fully discrete schemes, let us here
point out just one aspect, namely that the $H^1$ norm of the fully
discrete approximation is not entirely nonincreasing but can grow
slightly with a growth factor that, however, tends to zero as $\Dx\to
0$.  Compared to semi-implcit case \cite{Coclite:2006kx}, the proof is
notably more complicated and involves working with a version of the
scheme \eqref{eq:u-disc-fully} in which the quadratic terms have been
suitably truncated.

The paper is organized as follows: In Section \ref{sec:prelim} we
introduce some notation to improve the readability and recall a few
mathematical results relevant for the convergence analysis.  The
finite difference scheme and its convergence theorem are stated in
Section \ref{sec:diffscheme}. The convergence theorem is a consequence
of the results proved in Sections
\ref{sec:enest}-\ref{sec:strongconv}.  Finally, we present a numerical
example in Section \ref{sec:numex}.


Throughout this paper we use $C$ to denote a generic constant; the
actual value of $C$ may change from one line to the next in a
calcuation. We also use the notation that $a_{i} \lesssim b_i$ to mean
that $a_i \le Cb_i$ for some positive constant $C$ which is
independent of $i$.

\section{Preliminaries}\label{sec:prelim}

In what follows,  $\Dx$ and $\Dt$ denote two small positive numbers.
Unless otherwise stated, the indices $j$ and $n$ will run over $\Z$ and 
$0,\ldots N$, respectively, where $N\Dt = T$ for 
a fixed final time $T>0$. For such indices we set $x_j=j\Dx$, 
$x_{j+1/2}=(j+1/2)\Dx$, $t^n=n\Dt$, and introduce the grid cells
$$
I_j=[x_{j-1/2},x_{j+1/2}),\quad
I^n=[t^n,t^{n+1}),\quad\text{and}\quad I^n_j=I_j\times I^n.
$$

The following notations will be used frequently:
$$
a\vee 0 =\max\seq{a,0}=\frac{a+\abs{a}}{2}, 
\qquad a\wedge 0
=\min\seq{a,0}=\frac{a-\abs{a}}{2}.
$$

For $n\in \seq{0,\dots,N}$, let $v^n=\seq{v_j^n}_{j\in \Z}$ denote 
an arbitrary sequence, where $n$ refers to ``time'' and $j$ to ``space".  
We will frequently employ the following finite difference operators:
\begin{equation*}
  \begin{gathered}
    D_+v^n_j:=\frac{v^n_{j+1}-v^n_j}{\Dx},\qquad
    D_-v^n_j:=\frac{v^n_{j}-v^n_{j-1}}{\Dx},\\
    Dv^n_j:=\frac{D_+v^n_j+D_-v^n_j}{2}
    =\frac{v^n_{j+1}-v^n_{j-1}}{2\Dx}\quad\text{and}\quad \DT v^n_j
    =\frac{v^{n+1}_j - v^n_j}{\Dt}.
  \end{gathered}
\end{equation*}
We also use the notations
\begin{align*}
  \norm{v^n}_{\ell^p}& :=\left(\Delta x\sum\limits_{j\in\Z}
    \abs{v^n_j}^p\right)^{\frac1p},\, 1\le p<\infty, \qquad
  \norm{v^n}_{\ell^\infty} := \sup_j \abs{v_j^n}, 
  \\ \norm{v^n}_{h^1}&:=\left(\Dx\sum\limits_{j\in\Z} 
    \left[ \left(v^n_j\right)^2+ (D_-v^n_j)^2\right]\right)^{\frac12}.
\end{align*}
Occasionally, we also use the ``space-time''  $\el{p}$ norms of 
$v=\seq{v^n}_{n=0}^N=\seq{v_j^n}_{j,n}$:
$$
\norm{v}_{\el{p}} := \left(\Dt\sum_{n=0}^N 
\norm{v^n}_{\el{p}}^p\right)^{1/p}.
$$
Note that if $v\in\ell^p$, $p<\infty$, then $\lim_{j\to\pm\infty}
v_j=0$.

Let $\seq{v_j}_{j\in \Z}$ and $\seq{w_j}_{j\in \Z}$ denote two 
arbitrary (spatial) sequences.  Suppose $\norm{\seq{v_j}_j}_{h^1}<\infty$.  
Then the following discrete Sobolev inequality holds:
\begin{equation}\label{eq:discrete-sobol}
        \norm{\seq{v_j}_j}_{\ell^\infty}
        \le \frac{1}{\sqrt{2}}\norm{\seq{v_j}_j}_{h^1}.
\end{equation}
The discrete product rule takes the form
\begin{equation}\label{eq:discreteLeibnitz}
        \Dpm\left(v_jw_j\right)=v_j\Dpm w_j + \Dpm v_j w_{j\pm 1}.
\end{equation}
Moreover, the discrete chain rule states
\begin{equation}\label{eq:taylor1}
        \Dpm  f(v_{j})= f'(v_{j})\Dpm v_j
        \pm \frac\Dx2 f''(\xi_j^\pm)(\Dpm v_j)^2,       \qquad f\in C^2,
\end{equation}
for some number $\xi_j^\pm$ between $v_{j\pm 1}$ and $v_j$. 

We continue to collect some handy results for later use, starting with 
a discrete Gronwall inequality.

\begin{lemma}\label{lem:discGronwall} 
Assume that $c^k\ge 0$ and $f^k\ge 0$ for 
all $k=0,\dots,N$, and that the sequence $\seq{u^n}_{n=0}^N$ satisfies 
the difference inequality
\begin{equation}\label{eq:diffineq}
        \DT u^n + f^n \le c^n u^n, \qquad n=0,\dots,N-1.
\end{equation}
If $u^n\ge 0$ for all $n=0,\dots,N$, then
$$
u^N + \exp\left(\Dt \sum_{n=0}^{N-1} c^n\right) \Dt \sum_{n=0}^{N-1}
\exp\left(-\Dt \sum_{k=0}^n c^k\right)f^n \le 
\exp\left(\Dt \sum_{n=0}^{N-1} c^n\right) u^0.
$$
\end{lemma}
\begin{proof}
Set $R^n=\exp\left(-\Dt\sum_{k=0}^{n-1} c^k\right)$. Then we have
\begin{align*}
        \DT R^n &= \frac{1}{\Dt}\left(R^{n+1}-R^n\right)\\
        &=\exp\left(-\Dt\sum_{k=0}^{n} c^k\right)\frac{1}{\Dt}
        \left(1-\exp\left(\Dt c^n\right)\right)\\
        &\le R^{n+1} \frac{1}{\Dt}\left(1-(1+c^n\Dt)\right) 
        = -c^n R^{n+1}.
\end{align*}
Hence, multiplying \eqref{eq:diffineq} with $R^{n+1}$ we arrive at
$$
\DT\left(R^n u^n\right)= \DT u^n R^{n+1} + \DT R^{n} u^n 
\le -f^nR^{n+1}.
$$
Multiplying this by $\Dt$ and summing over $n$, we 
see that the lemma holds.
\end{proof}

The next lemma contains estimates for the solution of a
discrete version of the differential equation  $P-\pxx P=f$.

\begin{lemma}\label{lem:plemma}
Let $\seq{f_j}_{j\in\Z}$ be a sequence in $\el1\cap\el2$, 
and denote by $\seq{P_j}_{j\in\Z}$ the solution to 
the difference equation
\begin{equation}\label{eq:pprob}
        P_j - \Dm\Dp P_j = f_j, \qquad j\in\Z.
\end{equation}
Introducing the notations
\begin{equation*}
        h=\left(1+2\frac{1-e^{-\kappa}}{(\Dx)^2}\right)^{-1},\quad
        \kappa=\ln\left(1+\frac{\Dx^2}{2}+\frac{\Dx}{2}\sqrt{4+\Dx^2}\right),
\end{equation*}
the solution $\seq{P_j}_{j\in\Z}$ takes the form
\begin{equation}\label{eq:Green-dr}
        P_j=h \sum\limits_{i\in \Z}e^{-\kappa \abs{j-i}} f_i, \qquad j\in\Z.
\end{equation}
Moreover, the following estimates hold:
\begin{align}
  \label{eq:pbound1}
  \norm{\seq{P_j}_j}_{\ell^\infty}, \norm{\seq{P_j}_j}_{\el1}
  & \le C\norm{\seq{f_j}}_{\el1},\\
  \label{eq:pbound2}
  \norm{\seq{D_+ P_j}_j}_{\ell^\infty},
  \norm{\seq{\Dp P_j}_j}_{\el1} & 
  \le C \norm{\seq{f_j}_j}_{\el1},\\
  \norm{\seq{P_j}_j}_{h^1}&\le
  C\norm{\seq{f_j}_j}_{\el2},
  \label{eq:pbound3}
\end{align}
where $C>0$ is a constant independent of $\Dx$.
\end{lemma}
\begin{proof}
To verify the solution formula \eqref{eq:Green-dr}, we define $p_i$ by
$$
p_i=ce^{-\sigma \abs{i}},
$$
for some constants $c$ and $\sigma$ yet to be found. 
We shall choose these so that 
$$
(I-\Dp\Dm)p_i
=
\begin{cases}
        1, &\text{if $i=0$,}\\
        0, &\text{otherwise.}
\end{cases}
$$
If we find that this holds with $\sigma=\kappa$ and $c=h$ then
\eqref{eq:Green-dr} holds. We observe that for $i\ne 0$
$$
\Dp\Dm p_i =c e^{-\sigma\abs{i}} 2 \frac{\cosh(\sigma)-1}{\Dx^2}.
$$
Hence $\sigma$ must satisfy
$$
\sigma=\cosh^{-1}\left(1+\frac{\Dx^2}{2}\right) = \kappa.
$$
For $i=0$ we find that
$$
p_0-\Dp\Dm p_0 = c\left(1-\frac{2}{\Dx^2}\left(e^{-\kappa} -
1\right)\right). 
$$
If this is to be equal $1$ then $c=h$.  

For later use, one should observe that
\begin{equation}\label{eq:hkO}
        h=\frac{\Dx}{2}+{\mathcal O}(\Dx^2),\quad
        \frac{\abs{e^\kappa-1}}{\Dx}=1+\mathcal{O}(\Dx),
        \quad \frac{\abs{e^{-\kappa}-1}}{\Dx}=1+\mathcal{O}(\Dx).
\end{equation}

For any $j\in\Z$, we have
$$
\abs{P_j} \lesssim \norm{\seq{f_j}}_{\el1}
$$
Furthermore,
\begin{align*}
  \norm{\seq{P_j}}_{\el1} \le 2h\sum\limits_{i} \bigg[\Delta x\sum\limits_{j}
  e^{-\kappa \abs{j-i}}\bigg] \abs{f_i} 
  \lesssim \norm{\seq{f_i}}_{\el1},
\end{align*}
Hence, we have proved \eqref{eq:pbound1}.

From \eqref{eq:Green-dr},
\begin{align*}
  D_+P_j & =\frac{P_{j+1}-P_{j}}{\Dx} \\ 
  & =
  h\sum\limits_{i}\frac{e^{-\kappa\abs{i-j-1}}-e^{-\kappa\abs{i-j}}}{\Dx}
  f_i \\  
  & = h\sum\limits_{i=j}^{\infty}\frac{e^{-\kappa(i-j-1)}-e^{-\kappa(i-j)}}{\Dx}
  f_i + h
  \sum\limits_{i=-\infty}^{j-1}\frac{e^{\kappa(i-j-1)}-e^{\kappa(i-j)}}{\Dx}
  f_i \\  
  & = h\sum\limits_{i=j}^{\infty}e^{-\kappa(i-j)}\frac{e^{\kappa}-1}{\Dx} f_i 
  +h\sum\limits_{i=-\infty}^{j-1}e^{\kappa(i-j)}\frac{e^{-\kappa}-1}{\Dx}f_i. 
\end{align*}
Using \eqref{eq:hkO} we acquire from this the following two estimates:
\begin{equation*}
  \abs{\Dp P_j} 
  \lesssim h \sum\limits_{i}e^{-\kappa\abs{i-j}} \abs{f_i}
  \lesssim \norm{\seq{f_i}}_{\el1}
\end{equation*}
and
\begin{equation*}
  \norm{\seq{D_+P_j}_j}_{\el1} \lesssim h\Dx
  \sum\limits_{j,i}e^{-\kappa\abs{i-j}} \abs{f_i}
  \lesssim\norm{\seq{f_i}}_{\el1},
\end{equation*}
Therefore
\eqref{eq:pbound2} holds.

It remains to prove \eqref{eq:pbound3}. To this end, we multiply the 
equation \eqref{eq:pprob} by $\Dx P_j$ and 
perform a summation by parts to discover
\begin{align*}
  \norm{\seq{P_j}_j}_{h^1}^2 & 
  = \Dx\sum_j P_j f_j
  \le \frac12  \norm{\seq{P_j}_j}_{\el2}^2+\frac12\norm{\seq{f_j}_j}_{\el2}^2,
\end{align*}
from which \eqref{eq:pbound3} follows.
\end{proof}

We shall routinely use some well-known results related to weak
convergence, which we collect in a lemma (for proofs, see, e.g.,
\cite{Feireisl:2004oe}).  Throughout the paper we use overbars to
denote weak limits.
\begin{lemma}\label{lem:prelim} 
Let $O$ be a bounded open subset of $\R^M$, with $M\ge 1$.
  
Let $\Set{v_n}_{n\ge 1}$ be a sequence of measurable functions on
$O$ for which
$$
\sup_{n\ge 1} \int_O \Phi(\abs{v_n(y)})\, dy<\infty,
$$
for some given continuous function $\Phi:[0,\infty)\to[0,\infty)$.  
Then along a subsequence as $n \to \infty$
$$
\text{$g(v_n)\weakto \overline{g(v)}$ in $L^1(O)$}
$$
for all continuous functions $g:\R\to \R$ satisfying
$$
\lim_{\abs{v}\to \infty} \frac{\abs{g(v)}}{\Phi(\abs{v})}=0.
$$
Let $g\colon \R\to (-\infty,\infty]$ be a lower semicontinuous 
convex function and $\Set{v_n}_{n\ge 1}$ a sequence of measurable 
functions on $O$, for which
$$
\textit{$v_n\weakto v$ in $L^1(O)$, $g(v_n)\in L^1(O)$ for each 
$n$, $g(v_n)\weakto \overline{g(v)}$ in $L^1(O)$}.
$$
Then
$$
\text{$g(v)\le \overline{g(v)}$ a.e.~on $O$.}
$$
Moreover, $g(v)\in L^1(O)$ and
$$
\int_O g(v)\dy \le \liminf_{n\to\infty} \int_O g(v_n)\dy.
$$
If, in addition, $g$ is strictly convex on an open interval
$(a,b)\subset \R$ and
$$
\text{$g(v)=\overline{g(v)}$ a.e.~on $O$},
$$
then, passing to a subsequence if necessary,
$$
\text{$v_n(y)\to v(y)$ for a.e.~$y\in \Set{y\in O\mid v(y)\in (a,b)}$.}
$$
\end{lemma}

Let $X$ be a Banach space and denote by $X^\star$ its dual.  The space
$X^\star$ equipped with the weak-$\star$ topology is denoted by
$X^\star_{\mathrm{weak}}$, while $X$ equipped with the weak topology
is denoted by $X_{\mathrm{weak}}$. By the Banach-Alaoglu theorem, a
bounded ball in $X^\star$ is $\sigma(X^\star,X)$-compact.  If $X$
separable, then the weak-$\star$ topology is metrizable on bounded
sets in $X^\star$, and thus one can consider the metric space
$C\left([0,T];X^\star_{\mathrm{weak}}\right)$ of functions $v:[0,T]\to
X^\star$ that are continuous with respect to the weak topology. We
have $v_n\to v$ in $C\left([0,T];X^\star_{\mathrm{weak}}\right)$ if
$\langle v_n(t),\phi \rangle_{X^\star,X}\to \langle v(t),\phi
\rangle_{X^\star,X}$ uniformly with respect to $t$, for any $\phi\in
X$. The following lemma is a consequence of the Arzel\`a-Ascoli
theorem:

\begin{lemma}\label{lem:timecompactness}
Let $X$ be a separable Banach space, and suppose $v_n\colon [0,T]\to
X^\star$, $n=1,2,\dots$, is a sequence of measurable functions such
that
$$
\norm{v_n}_{L^\infty([0,T];X^\star)}\le C,
$$
for some constant $C$ independent of $n$. Suppose the sequence
$$
[0,T]\ni t\mapsto \langle v_n(t),\Phi \rangle_{X^\star,X}, 
\quad n=1,2,\dots,
$$
is equi-continuous for every $\Phi$ that belongs to a dense
subset of $X$.  Then $v_n$ belongs to
$C\left([0,T];X^\star_{\mathrm{weak}}\right)$ for every
$n=1,2,\dots$, and there exists a function $v\in 
C\left([0,T];X^\star_{\mathrm{weak}}\right)$ such that along a 
subsequence as $n\to \infty$
$$
\text{$v_n\to v$ in $C\left([0,T];X^\star_{\mathrm{weak}}\right)$}.
$$
\end{lemma}

\section{Explicit scheme and main result}\label{sec:diffscheme}
In this section we present the fully discrete (explicit) difference scheme for 
the Cammassa-Holm equation \eqref{eq:uprime}, which generates 
sequences $\seq{\unjph}$ and $\seq{\Pnj}$ 
for $(n,j)\in \seq{0,\dots,N}\times \Z$. We let $\seq{\unjph}$ solve the 
explicit difference equation
\begin{equation}\label{eq:u-disc}
  \DT\unjph + \left(\unjph\vee 0 \right)\Dm\unjph 
  + \left(\unjph\wedge 0 \right)\Dp\unjph +\Dp\Pnj=0,
\end{equation}
where the initial values are specified as follows:
\begin{equation}\label{eq:u-disc-data}
  \ujph^0=
  u_0(x_{j+1/2}),
\end{equation}
Given $\seq{\unjph}$, we determine $\seq{\Pnj}$ by solving
\begin{equation}\label{eq:P-disc}
  -D_-D_+\Pnj+\Pnj =\left(\unjph\vee 0 \right)^2 +\left(\unjmh\wedge 0\right)^2
  + \frac12 \left(\Dm\unjph\right)^2,
\end{equation}
which is a linear system of equations that can be solved 
as outlined in Lemma \ref{lem:plemma}.

Next, let us derive the difference scheme satisfied by
\begin{equation}\label{eq:qjdef}
        \qnj=\Dm \unjph.
\end{equation}
This will be done by applying the difference operator $\Dm$ to 
the $u$-equation \eqref{eq:u-disc}.  To this end, we apply the 
discrete product rule to find
$$
\Dm \left[\left(\unjph\vee 0 \right)\Dm\unjph\right] =
\left(\unjmh\vee 0 \right)\Dm \qj +\Dm\left(\unjph\vee 0 \right)\qnj
$$
and
$$
\Dm \left[\left(\unjph\wedge 0 \right)D_{+}\unjph\right] =
\left(\unjph\wedge 0 \right)D_{+}\qj +\Dm\left(\unjph\wedge 0 \right)\qj,
$$
so that
\begin{equation}\label{eq:q-der-I}
        \begin{split}
                & \Dm \left[\left(\unjph\vee 0 \right)\Dm\unjph 
                +\left(\unjph\wedge 0 \right)D_{+}\unjph\right] 
                \\ & \qquad \qquad 
                = \left(\unjmh\vee 0 \right)\Dm \qnj 
                + \left(\unjph\wedge 0 \right)D_{+}\qnj 
                + \left({\qnj}\right)^2.
        \end{split}
\end{equation}
The $P$-equation \eqref{eq:P-disc} rephrased in terms of $q$ reads
\begin{equation}\label{eq:q-der-II}
        -\Dm\Dp\Pnj+\Pnj
        = \left(\unjph\vee 0 \right)^2 
        +\left(\unjmh \wedge 0\right)^2 + \frac12 \left({\qnj}\right)^2.
\end{equation}
Employing \eqref{eq:q-der-I} and \eqref{eq:q-der-II} when applying
$\Dm$ to the $u$-equation in \eqref{eq:u-disc} yields
\begin{equation}\label{eq:q-disc}
        \begin{split}
                \DT\qnj & + \left(\unjmh\vee 0 \right)\Dm \qnj 
                  + \left(\unjph\wedge 0 \right)\Dp\qnj 
                \\ & + \frac{\left({\qnj}\right)^2}{2} +\Pnj 
                - \left(\unjph\vee 0\right)^2 
                - \left(\unjmh\wedge 0\right)^2 =0.
        \end{split}
\end{equation}

Regarding the initial values, in view of  \eqref{eq:qjdef} and \eqref{eq:u-disc-data}, 
we observe that
\begin{equation}\label{eq:init-q}
        \qj^0=\frac{1}{\Dx}\int_{I_j} \px u_0(x)\dx, \qquad j\in\Z.
\end{equation}

Since the variable $q=\px u$ can be discontinuous, \eqref{eq:q-disc} represents a 
natural upwind discretization of the 
equation for $q$, $\pt q + u \px q + \frac{q^2}{2}-u^2 + P=0$.

The main result of this paper is the convergence of the scheme to a
dissipative weak solution of \eqref{eq:u}-\eqref{ass:init}, which is
defined in the following sense \cite{Xin:2000qf,Xin:2002ys}:

\begin{definition}\label{def:sol}
Fix a final time $T>0$. We call a 
function $u\colon [0,T]\times\R\rightarrow \R$ a weak solution 
of the Cauchy problem for \eqref{eq:u}-\eqref{ass:init} on $\cDo$ if
\begin{Definitions}
\item $u\in C([0,T]\times\R)\cap L^\infty(0,T;H^1(\R))$;\label{def:solh1}
  
\item For all $s$ and $t$ in $[0,T]$, with $s\le t$, we have
  $\norm{u(t,\cdot)}_{H^1(\R)}\le  
  \norm{u(s,\cdot)}_{H^1(\R)}$. \label{def:solH1} 
                
\item $u$ satisfies \eqref{eq:uprime} in the sense of distributions on
  $\Do$;\label{def:soldistri} 
        
\item $u(0,x)= u_0(x)$ for every $x\in \R$;\label{def:initdata}
         
\item If, in addition, there exists a positive constant $K$ 
  such that
  \begin{equation*}
    u_x(t,x)\le \frac{2}{t}+ K \norm{u_0}_{H^1(\R)}^2, 
    \qquad  (t,x)\in (0,T]\times\R,
  \end{equation*}
  then we call $u$ a dissipative weak solution of the Cauchy problem
  \eqref{eq:u}-\eqref{ass:init}. \label{def:Oleinik}
\end{Definitions}
\end{definition}

In addition to $\px u\in L^\infty(0,T;L^2(\R))$, cf.~\ref{def:solH1},
the dissipative weak solutions $u$ that we construct in this paper
will possess an improved integrability property, namely $\px u \in
L^p_{\loc}(\Do)$ for $p<3$, i.e.,
$$
\int_0^T \int_a^b \abs{\px u}^{p}
\dx\dt\le C(a,b,T,p), 
\quad \text{$\forall a,b\in \R$, $a<b$.}
$$

To state our main convergence result and also for later use, we need to introduce 
some functions (interpolations of the difference approximations) 
that are defined at all points $(t,x)$ in the domain. 
We begin by defining the functions
$$
\left. \begin{aligned}
  q_j(t)&=\qnj + \left(t-t^n\right)\DT\qnj, \\
  u_{j+1/2}(t)&=\unjph+\left(t-t^n\right)\DT\unjph
\end{aligned}\right\} 
\quad\text{for $t\in I^n$.}
$$
With the aid of these we define
\begin{equation}\label{eq:qdxdef}
  \qdx(t,x)=q_j(t), \qquad (t,x)\in I_j^n
\end{equation}
and
\begin{equation}\label{eq:udxdef}
        \udx(t,x)=u_{j-1/2}(t)+\left(x-x_{j-1/2}\right)q_j(t),
        \qquad\text{for $(t,x)\in I_j^n$,}
\end{equation}
for $j\in \Z$, $n=0,\dots,N-1$. Note that $t\mapsto \udx(t,x)$ is a continuous function,
since $u_{j-1/2}(t)$ and $q_j(t)$ are continuous. 
Regarding the continuity in $x$ we have that
\begin{align*}
        \lim_{x\uparrow x_{j+1/2}}\udx(t,x) 
        &= u_{j-1/2}(t)+\left(x_{x+1/2}-x_{j-1/2}\right)q_j(t)\\
        &=\unjmh + \left(t-t^n\right)\DT\unjmh 
        +\Dx\left(\qnj+\left(t-t^n\right)\DT\qnj\right)\\
        &=\unjmh + \left(t-t^n\right)\DT\unjmh + \left(\unjph-\unjmh\right)
        \\&\qquad\quad +\left(t-t^n\right)\DT\left(\unjph-\unjmh\right) =u_{j+1/2}(t),
\end{align*}
and therefore $\udx$ is continuous, and furthermore
$\partial_x\udx = \qdx$ almost everywhere. Observe also that, due 
to \eqref{eq:init-q}, there holds 
$\qdx(0,x)\to \px u_0$ in $L^2(\R)$ as $\Dx\to 0$. 
Similarly to $\udx$, we define a function $\Pdx$ by
bilinear interpolation. First, let
$$
P_j(t)=P^n_j + (t-t^n)\DT\Pnj, \qquad t\in I^n,
$$
and then define
\begin{equation}\label{eq:Pdxdef}
        \Pdx(t,x)=P_j(t)+(x-x_j)\Dp P_j(t),\qquad (t,x)\in I_j^n,
\end{equation}
for $j\in \Z$, $n=0,\dots,N-1$.

We are now in a position to state our main result.
\begin{theorem}\label{th:main}
Suppose \eqref{ass:init} holds. Let $\seq{\udx}_{\Dx>0}$ 
be a sequence defined by \eqref{eq:udxdef} 
and \eqref{eq:u-disc}-\eqref{eq:qjdef}. 
Then, along a subsequence as $\Dx\dto 0$,
\begin{equation*}
        \text{$\udx \to u $ in $H_\loc^1(\Do)$},
\end{equation*}
where $u$ is a dissipative weak solution of the 
Cauchy problem \eqref{eq:u}-\eqref{ass:init}.
\end{theorem}

This theorem is a consequence of the results stated and proved in
Sections \ref {sec:enest}-\ref{sec:strongconv}. 

\section{Total energy estimate and some consequences}
\label{sec:enest}

The purpose of this section is establish a discrete total energy 
estimate for the difference scheme \eqref{eq:u-disc}-\eqref{eq:qjdef}.

\begin{lemma}\label{lem:enest} 
Assume that $\Dx$ and $\Dt$ are related through the CFL type condition
\begin{equation}\label{eq:cfl}
  \Dt<\frac{\log(1+\Dx^\theta)\,\Dx^2}
  {C\norm{u_0}_{H^1(\R)}^2\left(1+\Dx^2\right)},
\end{equation}
where $C$ is a constant (to be detailed in the proof of the lemma) that 
is independent of $\Dx$ and $u_0$ and $\theta>0$. Then, for any $N_0\in
\seq{0,\dots,N}$, and for all sufficiently small $\Dx$,
\begin{equation}\label{eq:enest}
  \norm{u^{N_0}}_{h^1}^2 
  + \Dx^2\Dt\sum_{n=0}^{N_0-1}\sum_{j\in \Z}
  \abs{\unjph}\left(\Dm\Dp \unjph\right)^2 
  \le e^{t^N\Dx^\theta} \norm{u^0}_{h^1}^2. 
\end{equation}
\end{lemma}
\begin{proof} 
For the proof of \eqref{eq:enest}, we shall need to introduce an auxiliary 
difference scheme. To this end, we start by defining the cut-off function
$$
f^M(u)=
\begin{cases}
  -M, \; & u<M,\\
  u, \; & u\in[-M,M],\\
  +M \;& u>M,
\end{cases}
$$
where $M>0$ is a fixed constant (to be determined later on).  
Now, let $\seq{\tunjph}$ and $\seq{\tPnj}$ solve the following 
system of difference equations: 
\begin{equation}\label{eq:tudef}
  \begin{aligned}
    \DT\tunjph &+ \left(\fm\left(\tunjph\right)\vee 0\right)\Dm\tunjph
    \\ & + \left(\fm\left(\tunjph\right)\wedge 0\right)\Dp\tunjph 
    +\Dp\tPnj = 0,
  \end{aligned}
\end{equation}
for $n=0,\dots,N-1$ and $j\in \Z$, and
\begin{equation*}
  \begin{aligned}
    -\Dm\Dp \tPnj + \tPnj &= \left(\fm\left(\tunjph\right)\vee 0\right)
    \left(\tunjph\vee 0\right)
    \\ &\qquad + \left(\fm\left(\tunjmh\right)\wedge 0\right)
    \left(\tunjmh\wedge 0\right)\\
    &\qquad + \frac{1}{2} \Dm\left(\fm\left(\tunjph\right)\right)
    \Dm\tunjph,
  \end{aligned}
\end{equation*}
for $n=0,\dots,N$ and $j\in \Z$. Regarding the initial data, we set 
$\tilde{u}_j^0=u_j^0$ for $j\in \Z$.

If we define $\tqnj:=\Dm\tunjph$, then it is straightforward to see that
$\seq{\tqnj}$ satisfies the difference equation
\begin{equation}\label{eq:tqdef}
  \begin{split}
    \DT\tqnj & + \left(\fm\left(\tunjmh\right)\vee 0\right)\Dm\tqnj
    +\left(\fm\left(\tunjph\right)\wedge 0\right)\Dp\tqnj 
    \\ & + \frac{1}{2}\Dm\left(\fm\left(\tunjph\right)\right)\tqnj
    \\ & - \left(\fm\left(\tunjph\right)\vee 0\right) \left(\tunjph\vee 0\right)
    \\ & - \left(\fm\left(\tunjmh\right)\wedge 0\right)
    \left(\tunjmh\wedge 0\right) + \tPnj = 0.
  \end{split}
\end{equation}

Multiplying \eqref{eq:tudef} by $\tunjph$ we find that
\begin{equation}
  \label{eq:uu-disc}
  \begin{split}
    & \tunjph\DT\tunjph + \left(\fm\left(\tunjph\right)\vee 0\right)
    \left(\tunjph\vee 0 \right)\tqnj
    \\ & \quad + \left(\fm\left(\tunjph\right)\wedge 0 \right)
    \left(\tunjph\wedge 0 \right)\tqnjpe +\Dp\tPnj \tunjph=0,
  \end{split}             
\end{equation}
while multiplying \eqref{eq:tqdef} by $\tqnj$ gives us
\begin{equation}
  \label{eq:qq-disc}
  \begin{split}
    \tqnj \DT\tqnj & + \left(\fm\left(\tunjmh\right)\vee 0 \right)(\Dm \tqnj) \tqnj 
    + \left(\fm\left(\tunjph\right)\wedge 0 \right)(\Dp\tqnj) \tqnj
    \\ & + \frac{1}{2}\Dm\left(\fm(\tunjph)\right) \left(\tqnj\right)^2 
    \\ & -\left(\fm\left(\tunjph\right)\vee 0\right) \left(\tunjph\vee 0\right) \tqnj 
    \\ & - \left(\fm\left(\tunjmh\right)\wedge 0\right) \left(\tunjmh\wedge 0\right) \tqnj 
    + \tPnj\tqnj=0.
  \end{split}
\end{equation}
Adding \eqref{eq:uu-disc} and \eqref{eq:qq-disc}, multiplying the
result with $\Dx$, and summing over $j$ yields
\begin{equation*}
  \Dx\sum_{j\in\Z}\left(\tunjph \DT\tunjph+\tqnj \DT \tqnj\right)
  + \text{I} + \text{II} + \text{III}=0,
\end{equation*}
where
\begin{align*}
  \text{I} &= \Dx \sum_{j\in\Z}\left(\fm\left(\tunjmh\right)\vee 0 \right)(\Dm \tqnj) \tqnj 
  \\ & \qquad+ \Dx \sum_{j\in\Z}
  \left(\fm\left(\tunjph\right)\wedge 0 \right) (\Dp\tqnj) \tqnj 
  +\frac{\Dx}{2} \sum_{j\in\Z} \Dm\left(\fm\left(\ujph\right)\right)
  \left(\tqnj\right)^2,
  \\ \text{II} &= \Dx \sum_{j\in\Z} \left(\fm\left(\tunjph\right)\vee 0 \right) 
  \left(\tunjph\vee 0 \right)\tqnj
  \\ &\qquad
  + \Dx \sum_{j\in\Z} \left(\fm\left(\tunjph\right)\wedge 0\right) 
  \left(\tunjph\wedge 0 \right)\tqnjpe 
  \\ & \qquad 
  -\Dx \sum_{j\in\Z} \left(\fm\left(\tunjph\right)\vee 0\right)
  \left(\tunjph\vee 0\right) \tqnj \\ & \qquad- \Dx\sum_{j\in\Z}
  \left(\fm\left(\tunjmh\right)\wedge 0\right) \left(\tunjmh\wedge 0\right) 
  \tqnj\equiv 0 \quad \text{(by shifting indices)},
  \\ \text{III} &= \Dx \sum_{j\in\Z} \Dp\tPnj \tunjph + \Dx\sum_{j\in\Z} 
  \tPnj \tqnj \equiv 0 \quad \text{(by summation by parts, cf.~\eqref{eq:qjdef}).}
\end{align*}
Let us now deal with term I. The discrete 
chain rule \eqref {eq:taylor1} tells us that 
$$
(\Dpm \tqnj) \tqnj = \Dpm\left(\frac{\left(\tqnj\right)^2}{2}\right)
\mp \frac{\Dx}{2} (\Dpm \tqnj)^2.
$$
Hence
\begin{align*}
  \text{I} &= \Dx \sum_{j\in\Z} \left(\fm\left(\tunjmh\right)\vee 0 \right) 
  \left[\Dm\left(\frac{(\tqnj)^2}{2}\right)+\frac{\Dx}{2} (\Dm \tqnj)^2\right] 
  \\ &\qquad + \Dx \sum_{j\in\Z}
  \left(\fm\left(\tunjph\right)\wedge 0 \right) 
  \left[ \Dp\left(\frac{(\tqnj)^2}{2}\right) - \frac{\Dx}{2} (\Dp \tqnj)^2\right] \\ & 
  \qquad + \frac{\Dx}{2} \sum_{j\in\Z}
  \Dm\left(\fm\left(\tunjph\right)\right) (\tqnj)^2 
  =\text{I}_{1} + \text{I}_{2},
\end{align*}
where
\begin{align*}
  \text{I}_{1} &= \Dx \sum_{j\in\Z} \left(\fm\left(\tunjmh\right)\vee 0 \right) 
  \Dm\left(\frac{\left(\tqnj\right)^2}{2}\right) \\ 
  &\qquad + \Dx \sum_{j\in\Z} \left(\fm\left(\tunjph\right)\wedge 0 \right) 
  \Dp\left(\frac{\left(\tqnj\right)^2}{2}\right) 
  \\ & \qquad +\frac{\Dx}{2} \sum_{j\in\Z} \Dm\left(\fm\left(\tunjph\right)\right)(\tqnj)^2, 
  \\ \text{I}_{2} &= \frac{\Dx^2}{2} \sum_{j\in\Z}
  \Bigl[\left(\fm\left(\tunjmh\right)\vee 0 \right)(\Dm\tqnj)^2 
  - \left(\fm\left(\tunjph\right)\wedge 0 \right)(\Dp \tqnj)^2\Bigr] 
  \\ & =\frac{\Dx^2}{2} \sum_{j\in\Z}
  \Bigl[\left(\fm\left(\tunjph\right)\vee 0 \right)(\Dp\tqnj)^2
  -\left(\fm\left(\tunjph\right)\wedge 0 \right)(\Dp \tqnj)^2\Bigr]\\
  &=\frac{\Dx^2}{2} \sum_{j\in\Z}
  \abs{\fm\left(\tunjph\right)}\left(\Dp \tqnj\right)^2\ge 0.
\end{align*}
To handle the $\text{I}_{1}$-term, we use the discrete 
product rule \eqref{eq:discreteLeibnitz}:
\begin{align*}
  & \Dm\left[\left(\fm\left(\tunjph\right)\vee 0 \right)
    \frac{\left(\tqnj\right)^2}{2}\right] \\ &\quad =
  \left(\fm\left(\tunjmh\right)\vee 0 \right)
  \Dm\left(\frac{\left(\tqnj\right)^2}{2}\right)
  +\Dm\left(\fm\left(\tunjph\right)\vee 0 \right)
  \frac{\left(\tqnj\right)^2}{2}, \\
  & \Dp\left[\left(\fm\left(\tunjmh\right)\wedge 0 \right)
    \frac{(\tqnj)^2}{2}\right] \\
  &\quad= \left(\fm\left(\tunjph\right)\wedge 0 \right)
  \Dp\left(\frac{\left(\tqnj\right)^2}{2}\right)
  +\Dp\left(\fm\left(\tunjmh\right)\wedge
    0\right)\frac{\left(\tqnj\right)^2}{2} \\ &\quad
  =\left(\fm\left(\tunjph\right)\wedge 0 \right)
  \Dp\left(\frac{\left(\tqnj\right)^2}{2}\right)
  +\Dm\left(\fm\left(\tunjph\right)\wedge 0\right)
  \frac{\left(\tqnj\right)^2}{2}.
\end{align*}
Using this we find that
\begin{align*}
  \text{I}_{1} &= \Dx \sum_{j\in\Z}
  \Dm\left[\left(\fm\left(\tunjph\right)\vee 0 \right)
    \frac{\left(\tqnj\right)^2}{2}\right] 
  \\ 
  & \qquad - \Dx
  \sum_{j\in\Z} \Dm\left(\fm\left(\tunjph\right)\vee 0 \right)
  \frac{\left(\tqnj\right)^2}{2} 
  \\ 
  &\qquad +
  \Dx\sum_{j\in\Z}\Dp\left[\left(\fm\left(\tunjmh\right)\wedge 0
    \right) \frac{\left(\tqnj\right)^2}{2}\right] 
  \\ 
  & \qquad -\Dx
  \sum_{j\in\Z} \Dm\left(\fm\left(\tunjph\right)\wedge 0 \right)
  \frac{\left(\tqnj\right)^2}{2} 
  \\
  &\qquad + \frac{\Dx}{2}
  \sum_{j\in\Z} \Dm\left(\fm\left(\tunjph\right)\right)
  \left(\tqnj\right)^2 
  \\ 
  &=-\Dx \sum_{j\in\Z} \Dm
  \left(\fm\left(\tunjph\right)\right) \frac{\left(\tqnj\right)^2}{2}
  \\ 
  & \qquad \qquad + \frac{\Dx}{2} \sum_{j\in\Z}
  \Dm\left(\fm\left(\tunjph\right)\right) \left(\tqnj\right)^2=0.
\end{align*}
Summarizing our findings so far:
\begin{equation}
  \label{eq:uu+qq-disc-end1}
  \begin{split}
    &\Dx\sum_{j\in\Z}\left(\tunjph\DT\tunjph+\tqnj \DT \tqnj\right)
    \\ & \qquad +\frac{\Dx^2}{2} \sum_{j\in\Z}\abs{\fm\left(\tunjph\right)}
    \left(\Dp \tqnj\right)^2=0.
  \end{split}
\end{equation}
Next, by \eqref{eq:taylor1},
\begin{equation}
  \label{eq:energy.0}
  \begin{split}
    &\Dx\sum_{j\in\Z}\left(\tunjph \DT\tunjph+\qj\DT\tqnj\right)
    \\ & \quad 
    =
    \DT\left[\frac{\Dx}{2}\sum_{j\in\Z}\left((\tunjph)^2+(\tqnj)^2\right)\right] 
    \\ &\quad\quad\quad 
    - \frac12\Dt \Dx\sum_{j\in\Z}\left(\left(\DT\tunjph\right)^2
      +\left(\DT\tqnj\right)^2\right).
  \end{split}
\end{equation}
Hence, we must now estimate
$$
\Dt\Dx\sum_{j\in\Z}\left( \left(\DT
\tunjph\right)^2+\left(\DT\tqnj\right)^2\right).
$$
Using \eqref{eq:tudef}, \eqref{eq:tqdef}, and the basic inequality
$\left(\sum_{\ell=1}^l a_l\right)^2\le 2^{l-1} \sum_{\ell=1}^l
(a_\ell)^2$, which holds for any sequence $\seq{a_\ell}_{\ell=1}^l$ of
positive real numbers, there is a positive constant $c_1$ that does
not depend on $\Dx$ such that
\begin{align*}
  & \Dt\Dx \sum_{j\in\Z}\left( \left(\DT
      \tunjph\right)^2+\left(\DT\tqnj\right)^2\right) 
  \\ 
  & \qquad \le
  c_1 \Dt\Dx\sum\limits_{j\in\Z} \Biggl[
  \left(\fm\left(\tunjph\right)\vee 0 \right)^2\left(\tqnj\right)^2 
  \\ 
  &\quad \qquad \qquad \qquad\qquad
  +\left(\fm\left(\tunjph\right)\wedge 0 \right)^2
  \left(\tilde{q}_{j+1}^n\right)^2+\left(\Dp\tPnj\right)^2 
  \\ 
  &\quad
  \qquad \qquad\qquad \qquad + \left(\fm\left(\tunjmh\right)\vee 0
  \right)^2 \left(\Dm \tqnj\right)^2 
  \\ 
  &\quad \qquad \qquad\qquad
  \qquad + \left(\fm\left(\tunjph\right)\wedge 0 \right)^2
  \left(\Dp\tqnj\right)^2
  \\
  & \quad \qquad \qquad\qquad \qquad
  +\left(\Dm\left(\fm\left(\tunjph\right)\right)\right)^2
  \left(\tqnj\right)^2 
  \\ 
  & \quad \qquad \qquad\qquad \qquad +
  \left(\fm\left(\tunjph\right)\vee 0\right)^2 \left(\tunjph\vee
    0\right)^2 
  \\ & \quad \qquad \qquad\qquad \qquad
  +\left(\fm\left(\tunjmh\right)\wedge 0\right)^2 \left(\tunjmh\wedge
    0\right)^2 + \left(\tPnj\right)^2\Biggr] 
  \\ 
  & \qquad \le
  c_1\Dt\left(J_1+J_2+J_3\right),
\end{align*}
where
\begin{align*}
  J_1&=
  \Dx\sum\limits_{j\in\Z}\Biggl[\left(\fm\left(\tunjph\right)\vee
    0\right)^2 \left(\tqnj\right)^2 +
  \left(\fm\left(\tunjph\right)\wedge 0\right)^2
  \left(\tilde{q}_{j+1}^n\right)^2 \\& \qquad\qquad\qquad
  +\left(\fm\left(\tunjph\right)\vee 0\right)^2 \left(\tunjph\vee
    0\right)^2 \\& \qquad\qquad\qquad\qquad\qquad
  +\left(\fm\left(\tunjmh\right)\wedge0\right)^2 \left(\tunjmh\wedge
    0\right)^2 \Biggr], \\ J_2 &=\Dx\sum\limits_{j\in\Z}\Biggl[
  \left(\fm\left(\tunjmh\right)\vee 0\right)^2 \left(\Dm
    \tqnj\right)^2 \\ & \qquad\qquad\quad
  +\left(\fm\left(\ujph\right)\wedge 0\right)^2
  \left(\Dp\tqnj\right)^2
  +\left(\Dm\left(\fm\left(\tunjph\right)\right)\right)^2\left(\tqnj\right)^2\Biggr], 
  \\ J_3 & =\Dx\sum\limits_{j\in\Z}\Biggl[\left(\Dp\tPnj\right)^2
  +\left(\tPnj\right)^2\Biggr].
\end{align*}
Since $\abs{\fm(u)}\le M$, the following bounds hold:
\begin{equation}\label{eq:J_1}
  J_1 \le
  2M^2\Dx\sum\limits_{j\in\Z}\left[\left(\tunjph\right)^2+(\tqnj)^2\right] 
  = 2M^2\norm{\tilde{u}^n}_{h^1}
\end{equation}
and
\begin{equation}\label{eq:J_2}
  \begin{split}
    J_2&=\frac{\Dx}{\Dx^2} \sum\limits_{j\in\Z}
    \Biggl[\left(\fm\left(\tunjmh\right)\vee 0\right)^2\left(\Dm\tqnj\Dx\right)^2 
    \\ & \qquad \qquad \qquad \qquad 
    + \left(\fm\left(\tunjph\right)\wedge 0\right)^2 \left(\Dp\tqnj\Dx\right)^2
    \\ & \qquad \qquad \qquad\qquad\qquad 
    +\left(\Dm\left(\fm\left(\tunjph\right)\right)\Dx\right)^2
    \left(\tqnj\right)^2\Biggr]
    \\ & \le c_2M^2 \frac{\Dx}{\Dx^2}\sum\limits_{j\in\Z}(\qj)^2 
    \le c_2\frac{M^2}{\Dx^2}\norm{\tilde{u}^n}_{h^1},
  \end{split}
\end{equation}
for some constant $c_2>0$ independent of $\Dx$. 

To estimate $J_3$ we use Lemma \ref{lem:plemma}, specifically
\eqref{eq:pbound3}, which implies that
\begin{equation}\label{eq:J_3}
  \begin{split}
    J_3&\le c_3\Dx\sum\limits_{j\in\Z}
    \Biggl[\left(\fm\left(\tunjph\right)\vel 0\right)^2
    \left(\tunjph\vel 0\right)^2  
    \\ & \qquad \qquad \qquad \qquad 
    +\left(\fm\left(\tunjmh\right)\et 0\right)^2 \left(\tunjmh\et 0\right)^2 
    \\ & \qquad \qquad \qquad\qquad\qquad 
    +\left(\Dm\left(\fm\left( \tunjph\right)\right)\right)^2
    \left(\tqnj\right)^2\Biggr] 
    \\ &\le 2c_3M^2\Dx\sum\limits_{i\in\Z}\left(\tunjph\right)^2
    +c_3M^2\frac{\Dx}{\Dx^2}\sum\limits_{i\in\Z}
    \left(\tqnj\right)^2
    \\ &\le c_3M^2\left(1+\frac{1}{\Dx^2}\right)
    \norm{\tilde{u}^n}_{h^1},
  \end{split}
\end{equation}
for some constant $c_3>0$ independent of $\Dx$.

Blending \eqref{eq:J_1}, \eqref{eq:J_2}, and \eqref{eq:J_3} we 
derive the bound
\begin{equation}\label{eq:uu+qq-disc-end2}
  \Dt\Dx\sum_{j\in\Z}\left(\left(\DT\tunjph\right)^2
    +\left(\DT\tqnj\right)^2\right)
  \le CM^2\left(\Dt+\frac{\Dt}{\Dx^2}\right)\norm{\tilde{u}^n}_{h^1},
\end{equation}
where the constant $C$ is independent of $\Dx$. 

Combining \eqref{eq:uu+qq-disc-end1}, \eqref{eq:energy.0}, and 
\eqref{eq:uu+qq-disc-end2}, it follows that $\seq{\tunjph}$ obeys the
following discrete energy estimate:
\begin{equation*}
        \DT\norm{\tilde{u}^n}_{h^1}^2
        + \underbrace{\frac{\Dx^2}{2} \sum_{j\in\Z}\abs{\fm\left(\tunjph\right)}
        \left(\Dp \tqnj\right)^2}_{=:Z^n} \le
        \underbrace{CM^2\left(\Dt+\frac{\Dt}{\Dx^2}\right)}_{=:\omega}
        \norm{\tilde{u}^n}_{h^1}^2. 
\end{equation*}
By the discrete Gronwall inequality, cf.~Lemma~\ref{lem:discGronwall}, 
\begin{equation*}
  \norm{\tilde{u}^{N_0}}_{h^1}^2 + \sum_{n=0}^{N_0-1}
  e^{CM^2\omega (t^{N_0-1}-t^n)}Z^n
  \le e^{CM^2\omega t^{N_0-1}}\norm{u^0}_{h^1}^2.
\end{equation*}
Choosing $M=\norm{u_0}_{H^1(\R)}$ and recalling the CFL type condition
\eqref{eq:cfl}, we deduce 
\begin{align*}
  e^{CM^2\omega t^N}&=e^{C\norm{u_0}_{H^1(\R)}^2\Dt^2((1+\Dx^2)/\Dx^2)N}\\
  &\le (1+\Dx^\theta)^{N\Dt} \\
  &\le e^{t^N\Dx^\theta}
  \le 2 \quad \text{if}\quad 
  \Dx^\theta \le \frac{\log 2}{T}. 
\end{align*}
Therefore, in particular
$$
\norm{\tilde{u}^n}_{h^1}\le \sqrt{2}\norm{u^0}_{h^1}
\quad\text{for $n=0,\dots,N_0$.}
$$
By the discrete Sobolev inquality \eqref{eq:discrete-sobol}, we
find that
$$
\norm{\tilde{u}^n}_{\el\infty} \le
\frac{1}{\sqrt{2}}\norm{\tilde{u}^n}_{h^1} \le \norm{u^0}_{h^1} \le M
\quad\text{for $n=0,\dots,N_0$.}
$$
This means that $\tilde{u}^n$ will ``never notice'' $\fm$, since
$$
\fm\left(\tunjph\right) 
= \tunjph\quad\text{for $j\in \Z$ 
and $n=0,\dots,N_0$.}
$$
Therefore,
$$
\tunjph=\unjph\quad\text{and}\quad \tPnj=\Pnj 
\quad\text{for $j\in \Z$ and $n=0,\dots,N_0$.}
$$
Finally, \eqref{eq:enest} follows by noting that
$$
e^{CM^2\omega (t^{N_0-1}-t^n)}\ge 1 \quad\text{for $n=0,\dots,N_0-1$.}
$$
\end{proof}

We conclude this section by stating some immediate
consequences of \eqref{eq:enest}.

\begin{lemma}\label{lem:boundP} 
For $n=0,\dots,N$,
\begin{align*}
	\norm{P^n}_{\el\infty}, \norm{P^n}_{\el1} & \le C\norm{u_0}^2_{H^1(\R)},\\
	\norm{\Dp P^n}_{\el\infty}, \norm{\Dp P^n}_{\el1} 
	& \le C \norm{u_0}^2_{H^1(\R)},
\end{align*}
where $C>0$ is a constant independent of $\Dx$.
\end{lemma}
\begin{proof}
This follows immediately from 
Lemma~\ref{lem:plemma}, noting that in this case
$$
f_j=\left(\unjph\vee 0 \right)^2 +\left(\unjmh\wedge 0 \right)^2
+ \frac12 \left(\Dm\unjph\right)^2,
$$
and thus
$$
\norm{f_j}_{\el1}\le \norm{u^n}^2_{h^1}\le (1+\Dx)^2
\norm{u^0}_{H^1(\R)}^2.
$$
\end{proof}

\section{One-sided sup-norm estimate}\label{sec:Oleinik}

\begin{lemma}\label{lem:Oleinik} 
Assume that $\Dt$ satisfies the CFL type condition \eqref{eq:cfl} and 
that $\Dx$ is sufficiently small. For $n=0,\dots,N$ and $j\in\Z$, we 
then have
\begin{equation}\label{eq:oleinikbnd}
  \qnj \le \frac{2}{t^n} + C\norm{u^0}_{h^1},
\end{equation}
where $C>0$ is a finite constant. 
\end{lemma}
\begin{proof}
We can write the difference equation for $\seq{\qnj}$, see \eqref{eq:q-disc}, as
\begin{equation*}
  \begin{split}
    q^{n+1}_j &= \qnj\left(1-\lambda a - \lambda b\right)
    +\qnjm\lambda a +\qnjp\lambda b - \Dt\frac{{\qnj}^2}{2}
    \\ &\qquad + \Dt\left(\left(\unjmh\vee 0\right)^2
      +\left(\unjph\wedge 0\right)^2 - \Pnj\right),
  \end{split}
\end{equation*}
where
$$
a=\lambda \left(\unjmh\vee 0\right), 
\quad
b=-\lambda \left(\unjph\wedge 0\right), 
\quad \lambda=\Dt/\Dx.
$$ 
Now we have uniform bounds on $\norm{u^n}_{\el\infty}$ 
and $\norm{P^n}_{\el\infty}$ and thus 
\begin{equation}
  \label{eq:q-disc-est1}
  q^{n+1}_j \le \qnj\left(1-\lambda a - \lambda b\right)
  +\qnjm\lambda a +\qnjp\lambda b - \Dt\frac{{\qnj}^2}{2} + \Dt L,
\end{equation}
for some finite constant $L\lesssim \norm{u^0}_{h^1}^2$. 

Set $\bar{q}^n_j =\max\seq{\qnj,\qnjm,\qnjp}$. We claim that
\begin{equation}\label{eq:claimq}
  q^{n+1}_j\le \bar{q}^n_j 
  -\Dt\frac{\left(\bar{q}^n_j\right)^2}{2}+\Dt L, 
\end{equation}
if $\Dt$ is chosen sufficiently small. 

First we choose $\Dt$ so small that $\lambda(a+b)<1/2$. 
Then if $\bar{q}^n_j=\qnj$ the claim 
follows immediately from \eqref{eq:q-disc-est1}. 

Next, assume $\bar{q}^n_j=\qnjm$. Then note that
\begin{align*}
  \left(\qnj\right)^2&=\left(\qnjm\right)^2
  +\left(\qnj\right)^2-\left(\qnjm\right)^2 
  \\&=\left(\qnjm\right)^2 + \frac{\qnj+\qnjm}{2}\Dx\Dm \qnj
  \\&=\left(\qnjm\right)^2 + \frac{u_{j+1/2}^n-u^n_{j-3/2}}{2}\Dm\qnj.
\end{align*}
Since $\Dm\qnj<0$, we find that
$$
\left(\qnj\right)^2 \le \left(\qnjm\right)^2
-\norm{u^n}_{\el\infty} \Dm\qnj.
$$
Using this we can rephrase \eqref{eq:q-disc-est1} as
\begin{align*}
  q^{n+1}_j&\le 
  \qnj\left(1-\lambda\left(a+\norm{u^n}_{\el\infty}+b\right)\right) 
  +\qnjm\lambda\left(a+\norm{u^n}_{\el\infty}\right)
  +\qnjp\lambda b
  \\ & \qquad\quad 
  - \Dt\frac{\left(\qnjm\right)^2}{2} +\Dt L
  \\ & \le \bar{q}^n_j - \Dt\frac{\left(\bar{q}^n_j\right)^2}{2}
  +\Dt L=:F\left(\bar{q}^n_j\right),
\end{align*}
if $\lambda \norm{u^n}_{\el\infty}<\frac{1}{2}$. The proof 
of \eqref{eq:claimq} if $\bar{q}^n_j=\qnjp$ is similar.
  
Note that $F'(q)=1-\Dt q$, and thus $F$ is increasing for $q<1/\Dt$. 
Furthermore, by the CFL type condition \eqref{eq:cfl}, 
$\Dt=\mathcal{O}(\Dx^3)$, and by the bounds 
on $\norm{u^n}_{\el\infty}$,
$$
\qnj\le \abs{\qnj}\le \mathcal{O}\left(\frac{1}{\Dx}\right)=
\mathcal{O}\left(\frac{1}{\Dt^{1/3}}\right) \le \frac{1}{\Dt},
$$
for sufficiently small $\Dt$. Therefore, setting\footnote{This maximum exists 
since $q^n\in \ell^2$.} $M^n=\max_j \qnj$, from \eqref{eq:claimq} we get
$$
M^{n+1}\le F\left(M^n\right).
$$
Now set $Z^n=M^n-\sqrt{2L}$. Then
\begin{align*}
  Z^{n+1}&\le F\left(Z^n+\sqrt{2L}\right)-\sqrt{2L}\\
  &=Z^n+\sqrt{2L} -\frac{\Dt}{2}\left(Z^n+\sqrt{2L}\right)^2 + L\Dt
  -\sqrt{2L}\\
  &=Z^n\left(1-\Dt\sqrt{2L}\right)-\frac{\Dt}{2}\left(Z^n\right)^2\\
  &\le Z^n-\frac{\Dt}{2}\left(Z^n\right)^2.
\end{align*}
Now, clearly if $Z^n\le 0$, then $Z^{n+1}\le 0$. Hence, if $Z^0\le 0$,
then $Z^n\le 0$ for all $n>0$. If $Z^n>0$, by \cite[page
271]{Smoller:1994zn},
$$
Z^n\le \frac{2}{t^n+1/Z^0}\le \frac{2}{t^n}.
$$
This finishes the proof.
\end{proof}

\section{Higher integrability estimate}
\label{sec:highint}
We begin this section by deriving a ``renormalized form'' of the
finite difference scheme for $q_j$, so let $f$ be a nonlinear function
(renormalization) of appropriate regularity and growth. Multiplying
\eqref{eq:q-disc} by $f'(\qj)$ and using the discrete chain rule,
which in the present context reads
\begin{equation*}
  \begin{split}
    f'(\qnj) \Dpm \qnj & = \Dpm f(\qnj)\mp
    \frac{\Dx}{2}f''\left(q^n_{j\pm 1/2}\right) 
    \left(\Dpm \qnj\right)^2,
    \\ f'(\qnj) \DT \qnj & = \DT f\left(\qnj\right)-\frac{\Dt}{2}
    f''\left(q^{n+1/2}_j\right) \left(\DT\qnj\right)^2,
  \end{split}
\end{equation*}
where $q^n_{j\pm 1/2}$ is a number between $q_{j}$ and $q_{j\pm 1}$, and
$q^{n+1/2}_j$ is a number between $\qnj$ and $q^{n+1}_j$. 
Multiplying the scheme \eqref{eq:q-disc} with $f'(\qnj)$ we obtain
\begin{equation}\label{eq:fq-disc}
        \begin{split}
                \DT f\left(\qnj\right) & + \left(\unjmh\vee 0 \right)\Dm f(\qnj) 
                +\left(\unjph\wedge 0 \right)\Dp f(\qnj)
                \\ & +\frac{\left(\qnj\right)^2}{2}f'(\qnj)
                + \left[\Pnj-\left(\unjph\vee 0\right)^2
                -\left(\unjmh\wedge 0\right)^2\right]f'(\qnj)
                \\ & \qquad + I_{\Dx,f'',j}=\frac{\Dt}{2}
                f''\left(q^{n+1/2}_j\right)\left(\DT\qnj\right)^2,
        \end{split}
\end{equation}
where
\begin{align*}
        I_{\Dx,f'',j} & := \frac{\Dx}{2}\Biggl\{\left(\unjmh\vee 0 \right)
        f''(q^n_{j-1/2})(\Dm \qnj)^2 
        \\ & \qquad\qquad\quad - \left(\unjph\wedge 0 \right)
        f''(q^n_{j+1/2})(\Dp \qnj)^2\Biggr\}.
\end{align*}

Let us now write \eqref{eq:fq-disc} in divergence-form. To this end, observe 
that the discrete product rule \eqref{eq:discreteLeibnitz} implies 
the following relations:
\begin{align*}
        \Dm\left[ \left(\unjph\vee 0 \right) f(\qnj)\right] 
        &=\left(\unjmh\vee 0 \right) \Dm f(\qnj) + \Dm \left(\unjph\vee 0 \right) f(\qnj), 
        \\  \Dp\left[\left(\unjmh\wedge 0 \right) f(\qnj)\right]
        & = \left(\unjph\wedge 0 \right) \Dp f(\qnj) 
        + \Dp \left(\unjmh\wedge 0 \right) f(\qnj)
        \\ &= \left(\unjph\wedge 0 \right) \Dp f(\qnj)
        + \Dm \left(\unjph\wedge 0\right) f(\qnj),
\end{align*}
and therefore, using that $\qnj=\Dm\unjph$,
\begin{align*}
        & \left(\unjmh\vee 0 \right)\Dm f(\qnj) 
        + \left(\unjph\wedge 0\right)\Dp f(\qnj) 
        \\ & \qquad = \Dm\left[ \left(\unjph\vee 0 \right) f(\qnj)\right] 
        + \Dp\left[\left(\unjmh\wedge 0 \right)
        f(\qnj)\right] -\qnj f(\qnj).
\end{align*}
Hence, we end up with the following divergence-form variant of the
renormalized difference scheme \eqref{eq:fq-disc}:
\begin{equation} \label{eq:fq-disc-II}
        \begin{split}
                \DT f(\qnj)&+\Dm\left[ \left(\unjph\vee 0 \right)f(\qnj)\right] 
                + \Dp\left[\left(\unjmh\wedge 0 \right)f(\qnj)\right] 
                \\ & + \frac{\left(\qnj\right)^2}{2}f'(\qnj)-\qnj f(\qnj)
                \\ & \qquad +\left[ \Pnj - \left(\unjph\vee 0\right)^2 
                - \left(\unjmh\wedge 0\right)^2\right]f'(\qnj)
                \\ & \qquad\qquad +I_{\Dx,f'',j}=\frac{\Dt}{2}
                f''\left(q^{n+1/2}_j\right)\left(\DT\qnj\right)^2.
        \end{split}
\end{equation}

To ensure that the limit of $(\px \udx)^2$, cf.~ \eqref{eq:enest}, is nonsingular
(i.e., not a measure), we shall need the following higher integrability estimate:

\begin{lemma}\label{lem:highint} 
Let $\qdx$ be defined by \eqref{eq:qdxdef}, and assume that the 
CFL type condition \eqref{eq:cfl} holds. Then, for all finite 
numbers $a,b,\alpha$ with $a<b$ and $\alpha\in (0,\theta)$, 
\begin{equation}\label{eq:highint}
  \int_0^T\int_a^b \abs{\qdx}^{2+\alpha}\dx \dt\le C,
\end{equation}
for some constant $C=C(a,b,T,\alpha)$ that is independent of $\Dx$.
\end{lemma}
\begin{proof}
Define $\ee(q)=\sqrt{\eps^2+q^2}-\eps$. Note that 
$\ee(q) \approx \abs{q}$ for small $\eps$ and
$$
-1< \ee'(q)=\frac{q}{\sqrt{\eps^2+q^2}}<1\quad\text{and}\quad
0 < \ee''(q)=\frac{\eps^2}{(\eps^2+q^2)^{3/2}}\le \frac{1}{\eps}.
$$
In \eqref{eq:fq-disc-II}, we then specify
\begin{equation*}
        f(q)=\ee(q)\left(1+q^2\right)^{\alpha/2}.
\end{equation*}
One can easily check that 
\begin{align*}
        f'(q) &=\ee'(q)\left(1+q^2\right)^{\alpha/2}
        +\alpha \ee(q) q \left(1+q^2\right)^{\alpha/2-1}
        \\ f''(q)&=\ee''(q)\left(1+q^2\right)^{\alpha/2}
        + 2\alpha\ee'(q) q \left(1+q^2\right)^{\alpha/2-1}
        \\ & \qquad\quad +\alpha\ee(q)\left(\left(1+q^2\right)^{\alpha/2-1}
        +2\left(\frac{\alpha}{2}-1\right)
        \left(1+q^2\right)^{\alpha/2-2}\right),
\end{align*}
so that in particular $f''(q)\ge 0$ and 
\begin{equation}\label{eq:fppbnd}
        f''(q)=\ee''(q)\left(1+q^2\right)^{\alpha/2}
        + \text{bounded terms.}
\end{equation}

Next set
\begin{align*}
        H(q) & :=\frac{q^2}{2}f'(q) - qf(q)
        \\ & =\frac{q^2}{2}\left(\ee'(q)\left(1+q^2\right)^{\alpha/2}
        +\alpha\ee(q)\left(1+q^2\right)^{\alpha/2-1} q\right)
        -q\ee(q)\left(1+q^2\right)^{\alpha/2}
        \\& =q\ee(q)\left(1+q^2\right)^{\alpha/2}
        \left[\frac{q\ee'(q)}{2\ee(q)} + \frac{\alpha q^2}{2(1+q^2)}-1\right]
        \\ &=q\ee(q)\left(1+q^2\right)^{\alpha/2}\frac{1}{2}
        \left[\frac{\alpha q^2}{1+q^2}+\frac{\eps}{\sqrt{\eps^2+q^2}}-1\right]
        =:H_{\eps}(q)h_{\eps}(q),
\end{align*}
with
$$
H_{\eps}(q)=q\ee(q)\left(1+q^2\right)^{\alpha/2}, \qquad 
h_{\eps}(q)= \frac{1}{2}\left[\frac{\alpha q^2}{1+q^2}
+\frac{\eps}{\sqrt{\eps^2+q^2}}-1\right].
$$
Now note that
$$
\lim_{q\to-\infty}h_{\eps}(q)\le \lim_{q\to
  -\infty}h_1(q)=\frac{\alpha-1}{2}<0.
$$
Hence, for $\alpha,\eps<1$, we can find a constant $K>0$ such
that
\begin{equation}
  \label{eq:Kdef}
  h_{\eps}(q)<\frac{\alpha-1}{4}\quad \text{for all $q<-K$}.
\end{equation}

Let us continue by defining the sets
$$
\begin{gathered}
        \mathcal{N}^n=\seq{j\in\Z\;\bigm|\; \qnj<-K},\quad 
        \mathcal{C}^n=\seq{j\in\Z\;\bigm|\;-K\le \qnj\le 0},\\
    \text{and}\quad \mathcal{P}^n=\seq{j\in\Z\;\bigm|\; \qnj>0},
\end{gathered}
$$
where $K$ is defined in \eqref{eq:Kdef}. Moreover, let $0\le \chi(x)\le 1$ be a 
smooth cutoff function satisfying
$$
\chi(x)=
\begin{cases}
        0, &x<a-1, \\
        1, & x\in[a,b],\\
        0, & x>b+1.
 \end{cases}
$$
We multiply \eqref{eq:fq-disc-II} by $\chi(x_{j})\Dt\Dx$ 
and sum over $(n,j)\in \seq{0,\dots,N-1}\times \Z$ to get
\begin{align}
  \Dt\Dx\sum_{n}\sum_{j\in \mathcal{N}^n}
  h_{\eps}&\left(\qnj\right)H_{\eps}\left(\qnj\right)\chi(x_j)\notag 
  \\ 
  & \le -\Dt\Dx\sum_{n}\sum_{j\in \mathcal{C}^n}
  h_{\eps}\left(\qnj\right)H_{\eps}\left(\qnj\right)\chi\left(x_{j}\right) \label{eq:Cterm}
  \\ &\quad -\Dt\Dx\sum_{n}\sum_{j\in \mathcal{P}^n}
  h_{\eps}\left(\qnj\right)H_{\eps}\left(\qnj\right)\chi\left(x_{j}\right)
  \label{eq:Pterm}  
  \\ &  \quad 
  \begin{aligned}
    -\Dt\Dx\sum_{j,n}\Bigl[&\Dm\left(\left(\unjph\vee
        0\right)f\left(\qnj\right)\right)\chi(x_j)\\
    &\quad
    \Dp\left(\left(\unjmh\wedge
        0\right)f\left(\qnj\right)\right)\chi(x_j)\Bigr]
  \end{aligned} 
  \label{eq:diffterm}
  \\ 
  & \quad-\Dt\Dx\sum_{n,j}A^n_j\chi\left(x_{j}\right)
  \label{eq:Aterm}
  \\
  &\quad  +\Dx\sum_j
  \left(f\left(q^0_j\right)-f\left(q^N_j\right)\right) \chi(x_j)
  \label{eq:fterm} 
  \\ 
  &\quad +\Dx\Dt\sum_{n,j}\frac{\Dt}{2}f''\left(q^{n+1/2}_j\right)\chi(x_j)
  \left(\DT\qnj\right)^2,\label{eq:fppterm}
\end{align}
where 
\begin{equation}\label{eq:Anj}
        A^n_j=\Pnj-\left(\unjph\vee 0\right)^2
        -\left(\unjmh\wedge 0\right)^2.
\end{equation}
  
Now for $j\in\mathcal{N}^n$ we have that
$$
\frac{1-\alpha}{4}
\abs{\qnj}\ee\left(\qnj\right)\left(1+\left(\qnj\right)^2\right)^{\alpha/2}
\le h_\eps\left(\qnj\right)H_{\eps}\left(\qnj\right).
$$
Therefore
\begin{equation*}
        \begin{aligned}
                \frac{1-\alpha}{4}&\Dt\Dx \sum_{n}\sum_{j\in\mathcal{N}^n}
                \abs{\qnj}\ee\left(\qnj\right)
                \left(1+\left(\qnj\right)^2\right)^{\alpha/2}\chi(x_j)
                \\
                &\le \abs{\text{\eqref{eq:Cterm}}}
                +\abs{\text{\eqref{eq:Pterm}}} +\abs{\text{\eqref{eq:diffterm}}}
                +\abs{\text{\eqref{eq:Aterm}}} +\abs{\text{\eqref{eq:fterm}}}
                +\abs{\text{\eqref{eq:fppterm}}}.
        \end{aligned}
  \end{equation*}
We shall now find bounds on all the terms on the right hand side; in what follows, we 
let $C$ denote a generic constant independent of $\Dx$, $\eps$, and $\alpha$.
  
We start with \eqref{eq:fppterm}. By \eqref{eq:q-disc}
\begin{align*}
  \left(\DT \qnj\right)^2 & \le 
  C\Biggl[ \left(\left(\unjmh\vee 0\right)\Dm\qnj\right)^2
  + \left(\left(\unjph\wedge 0\right)\Dp\qnj\right)^2
  \\ & \qquad \qquad + \left(\Pnj\right)^2 
  +\left(\unjph\vee 0\right)^4 + \left(\unjmh\wedge 0\right)^4
  +\left(\qnj\right)^4\Biggr].
\end{align*}
To bound the ``integrals'' of these terms we must use the CFL type
condition \eqref{eq:cfl}, which implies that $\Dt=\mathcal{O}(\Dx^{2+\theta})$.
First,
\begin{align*}
  & \Dx \Dt \sum_j \left(\left(\unjmh\vee 0\right)\Dm\qnj\right)^2
  \\ & \quad \le C \Dx^{1+\theta}\sum_j \left(\Dx\Dm\qnj\right)^2
  \le C\Dx^{1+\theta} \sum_j \left(\qnj\right)^2 \le C \Dx^\theta.
\end{align*}
Therefore
\begin{align*}
  &\Dt\Dx\sum_{n,j} \Dt \Biggl[\left(\left(\unjmh\vee 0\right)\Dm\qnj\right)^2
  + \left(\left(\unjph\wedge 0\right)\Dp\qnj\right)^2 \Biggr]\chi(x_j)
  \\ & \qquad 
  \le C\Dx^\theta T \to 0 \quad \text{as $\Dx\to 0$.}  
\end{align*}
We also find that
\begin{align*}
  & \Dx\sum_j \Biggl[\left(\Pnj\right)^2+\left(\unjph\vee 0\right)^4
  +\left(\unjmh\wedge 0\right)^4\Biggr]
  \\ & \qquad \le C\Dx\sum_j \Biggl [ \abs{\Pnj}
  +\left(\unjph\vee 0\right)^2+\left(\unjmh\wedge 0\right)^2\Biggr]
  \\ & \qquad \le C\Dx \sum_j
  \Biggl[\abs{\Pnj}+\left(\unjph\right)^2\Biggr]\le C, 
\end{align*}
since $\udx$ and $P^n_j$ are uniformly bounded. Thus
$$
\Dt\Dx\sum_{n,j} \Dt \left[ \left(\Pnj\right)^2
+\left(\unjph\vee 0\right)^4 + \left(\unjmh\wedge 0\right)^4\right]
\chi(x_j)\to 0,
$$
as $\Dx \to 0$. Additionally,
\begin{align*}
  \Dx\sum_j \Dt \left(\qnj\right)^4 
  & \le C \Dx\sum_j \Dx^\theta \left(\Dx\qnj\right)^2 \left(\qnj\right)^2
  \\ & \le C\Dx^{1+\theta}\sum_j \left[\left(\unjph\right)^2
    +\left(\unjmh\right)^2\right]\left(\qnj\right)^2
  \\ & \le C\Dx^{1+\theta} \sum_j \left(\qnj\right)^2 
  \le C\Dx^\theta.
\end{align*}
Therefore
$$
\Dx\Dt\sum_j \Dt \left(\qnj\right)^4 \chi(x_j) \le CT\Dx^\theta \to 0 
\quad \text{as $\Dx\to 0$.} 
$$
Now we have established that
\begin{equation}\label{eq:dtsbound}
  \Dx\Dt\sum_{n,j}\Dt\left(\DT \qnj\right)^2\chi(x_j) 
  =\mathcal{O}\left(\Dx^\theta\right) 
  \quad \text{as $\Dx\to 0$.} 
\end{equation}  
Recalling \eqref{eq:fppbnd} this implies that
$$
\abs{\text{\eqref{eq:fppterm}}}\le C\Dx\Dt\sum_{n,j}
\frac{\Dt}{\eps}\left(1+\left(q^{n+1/2}_j\right)^2\right)^{\alpha/2}
\left(\DT\qnj\right)^2\chi(x_j) + \mathcal{O}\left(\Dx^\theta\right).
$$
When we established \eqref{eq:dtsbound} we always 
had a ``$\Dx^\theta$ to spare'', which we can use now. 
With $\beta=\theta-\alpha>0$, we get
\begin{align*}
  \Dt\Dx\sum_{n,j}
  &\frac{\Dt}{\eps}\left(1+\left(q^{n+1/2}_j\right)^2\right)^{\alpha/2}
  \left(\DT\qnj\right)^2\\ 
  & \le C \Dx\Dt\sum_{n,j} \frac{\Dx^{\beta}}{\eps}
  \Dx^{\alpha}\left(1+\left(q^{n+1/2}_j\right)^2\right)^{\alpha/2}
  \Dx^2\left(\DT\qnj\right)^2\\ 
  & \le C \Dx\Dt\sum_{n,j} \frac{\Dx^{\beta}}{\eps}
  \left(\Dx^2+\left(\Dx q^{n+1/2}_j\right)^2\right)^{\alpha/2}
  \Dx^2\left(\DT\qnj\right)^2 
  \\ & \le C\frac{\Dx^{\beta}}{\eps} \Dx\Dt\sum_{n,j}
  \Dx^2\left(\DT\qnj\right)^2
  \le C\frac{\Dx^{\beta}}{\eps},
\end{align*}
since
$$
\left(\Dx^2+\left(\Dx q^{n+1/2}_j\right)^2\right)^{\alpha/2} 
\le \left(\Dx^2+4\left(\norm{u^n}_{\el\infty}\right)^2\right)^{\alpha/2}\le C.
$$
Now choosing $\eps=\Dx^{\beta}$, we finally 
conclude that $\abs{\eqref{eq:fppterm}}$ is bounded.

Next we turn to \eqref{eq:Cterm}. For $-K\le q\le 0$, we 
have that $\abs{h_{\eps}(q)H_{\eps}(q)}\le C$, where $C$ is 
independent of $\eps$. Therefore
\begin{equation*}
        \abs{\text{\eqref{eq:Cterm}}}\le C T (b-a+2).
\end{equation*}

To estimate $\abs{\text{\eqref{eq:Pterm}}}$, observe that 
$$
\abs{\text{\eqref{eq:Pterm}}}
\le C\Dt\Dx\sum_{n}\sum_{j\in \mathcal{P}^n}
\left(1+\abs{\qnj}^{2+\alpha}\right)\chi\left(x_{j}\right)
$$
Now, equipped with \eqref{eq:enest} and  \eqref{eq:oleinikbnd}, it 
is possible to bound the right-hand side by a ``$\Dx$ independent" constant 
exactly as was done in \cite{Coclite:2006kx}.

Regarding \eqref{eq:diffterm}, observe that
\begin{align*}
        \biggl|\Dt & \Dx\sum\limits_{j,n} \Dm\left[\left(\unjph\vee 0\right)
        f\left(\qnj\right)\right]\chi(j\Dx)\biggr|
        \\ & =\bigg|\Dt\Dx\sum\limits_{j,n}\left(\unjph\vee 0\right)
        \Dp\chi(j\Dx)f\left(\qnj\right)\biggr| 
        \\ & \le C \Dt\Dx\sum\limits_{j,n} \abs{\unjph}\abs{\Dp\chi(j\Dx)}
        \left(1+\abs{\qnj}^{1+\alpha}\right) 
        \\ & \le C T \Biggl( \sup\limits_n \norm{u^n}_{\el\infty}
        \norm{\seq{\Dp\chi(x_j)}_j}_{\el1}
        \\ & \qquad\qquad\qquad \qquad
        +\sup\limits_n\norm{\seq{\abs{\qnj}^{1+\alpha}}_j}_{\el{\frac{2}{1+\alpha}}}
        \norm{\seq{\Dp\chi(j\Dx)}_j}_{\el{\frac{2}{1-\alpha}}}\Biggr)
        \\ & =c T \Biggl( \sup\limits_n\norm{u^n}_{\el\infty}
        \norm{\left\{\Dp\chi(x_j)\right\}_j}_{\el1}
        +\sup\limits_n\norm{q^n}^{1+\alpha}_{\el2}
        \norm{\seq{\Dp\chi(j\Dx)}_j}_{\el{\frac{2}{1-\alpha}}}\Biggr).
\end{align*}
Therefore, also $\abs{\eqref{eq:diffterm}}$ is bounded independently of $\Dx$.

Next we focus on \eqref{eq:Aterm}. Remembering 
that $\abs{f'(q)}\le C(1+\abs{q}^\alpha)$, we find
\begin{align*}
        \abs{\text{\eqref{eq:Aterm}}}
        & \le\Dt\Dx\sum\limits_{n,j}
        \left(\norm{P^n}_{\el\infty}+2\norm{u^n}_{\el\infty}^2\right)
        \abs{f'(\qnj)}\chi(x_j)
        \\ & \le C \Dt\Dx\sum\limits_{n,j}
        \left(\norm{P^n}_{\el\infty}+2\norm{u^n}_{\el\infty}^2\right)
        \left(1+\abs{\qnj}^{\alpha}\right)\chi(x_j)\\ 
        & \le C T\left(\sup\limits_n\norm{P^n}_{\el\infty}
        +2\sup\limits_n\norm{u^n}_{\el\infty}^2\right)
        \\ & \qquad\quad \times
        \left(\norm{\seq{\chi(x_j)}_j}_{\el1}
        +\sup\limits_n\norm{\seq{\abs{\qnj}^{\alpha}}_j}_{\el{\frac{2}{\alpha}}}
        \norm{\seq{\chi(j\Dx)}_j}_{\el{\frac{2}{2-\alpha}}}\right)
        \\ & \le C T\left(\sup\limits_n\norm{P^n}_{\el\infty}
        +2\sup\limits_n\norm{u^n}_{\el\infty}^2\right)
        \\ & \qquad \quad \times
        \left(\norm{\seq{\chi(j\Dx)}_j}_{\el1}
        +\sup\limits_n\norm{q^n}^{\alpha}_{\el2}
        \norm{\seq{\chi(j\Dx)}_j}_{\el{\frac{2}{2-\alpha}}}\right).
\end{align*}

Finally, keeping in mind that $f\ge 0$,  we treat \eqref{eq:fterm} as follows:
\begin{align*}
        \abs{\text{\eqref{eq:fterm}}}
        &\le \Dx\sum_j \left(f\left(q^N_j\right)+f\left(q^0_j\right)\right)\chi(x_j)
        \\ &\le C \Dx\sum_j \left(2+\abs{q^N_j}^\alpha
        +\abs{q^0_j}^\alpha\right)\chi(x_j)
        \\ & \le \norm{\seq{\chi(x_j)}_j}_{\el1}
        \\ & \quad 
        +\left(\norm{\seq{\abs{q^N_j}^{1+\alpha}}_j}_{\el{\frac{2}{1+\alpha}}}
        +\norm{\seq{\abs{q^0_j}^{1+\alpha}}_j}_{\el{\frac{2}{1+\alpha}}}\right)
        \norm{\seq{\chi(x_j)}_j}_{\el{\frac{2}{1-\alpha}}}
        \\ & = \norm{\seq{\chi(x_j)}_j}_{\el1} 
        + \left(\norm{q^N}_{\el2}^{1+\alpha}+ \norm{q^0}_{\el2}^{1+\alpha}\right)
        \norm{\seq{\chi(x_j)}_j}_{\el{\frac{2}{1-\alpha}}}.
\end{align*}
  
Summarizing, we have established 
\begin{equation*}
        \Dt\Dx \sum_n\sum_j \abs{\qnj}\ee\left(\qnj\right)
        \left(1+\left(\qnj\right)^2\right)^{\alpha/2}\le C.
\end{equation*}
The statement of the lemma follows by noting that
$$
\abs{q}^{2+\alpha}\le \abs{q}\ee(q)\left(1+q^2\right)^{\alpha/2}
+ \abs{q}^{1+\alpha},
$$
and using, in combination with \eqref{eq:enest}, the bound
$$
\Dt\Dx\sum_{j,n}\abs{\qnj}^{1+\alpha}\chi(x_j) 
\le CT\sup_n\norm{q^n}_{\el2}^{1+\alpha}
\norm{\seq{\chi(x_j)}_j}_{\el{\frac{2}{1-\alpha}}}.
$$
\end{proof}

\section{Basic convergence results}\label{sec:converg}

The purpose of this section is to present some straightforward 
consequences of the a priori estimates established in the foregoing 
sections. More precisely, we prove that the two sequences
$\seq{\udx}_{\Dx>0}$, cf.~\eqref{eq:udxdef}, and 
$\seq{\Pdx}_{\Dx>0}$, cf.~\eqref{eq:Pdxdef}, have strongly 
converging subsequences, starting with the former.

\begin{lemma}\label{lem:convlemma1}
There exists a limit function
$$
u\in L^\infty\left(0,T;H^1(\R)\right)\cap C(\cDo),
$$
such that along a subsequence as $\Dx\to 0$
\begin{align}
  \udx & \weakstar u \quad\text{in $L^\infty\left(0,T;H^1(\R)\right)$,}
  \label{eq:uweakstar} \\
  \udx &\to u \quad\text{uniformly in $[a,b]\times[0,T]$, for any $a<b$.}
  \label{eq:uuniform}
\end{align}
Additionally,
\begin{align}
  \label{eq:hdecrease}
  & \text{$t \mapsto \norm{u(t,\cdot)}_{H^1(\R)}$ is non-increasing, and}\\
  \label{eq:uinitialconv}
  &\lim_{t\to 0}u(t,x)=u_0(x), \qquad x\in \R.
\end{align}
\end{lemma}

\begin{proof}
First we note that for $t\in I^n$, 
\begin{align*}
  \int_{\R} (\udx(t,x))^2\dx & =\sum_j\! \int_{x_{j-1/2}}^{x_{j+1/2}}
  \!\!\! \left(\frac{1}{\Dx}(x_{j+1/2}-x)u_{j-1/2}(t)
    \! + \! (x-x_{j-1/2})u_{j+1/2}(t)\right)^2\dx \\ & 
  \le \frac{\Dx}{2}\sum_j \left( \left(u_{j-1/2}(t)\right)^2
    +\left(u_{j+1/2}(t)\right)^2\right) 
  \\ & \le \frac{1}{\Dt}\left(\left(t-t^n\right) \Dx\sum_j 
    \left(u^{n+1}_{j+1/2}\right)^2+ \left(t^{n+1}-t\right)\Dx\sum_j
    \left(\unjph\right)^2\right)
\end{align*}
and
\begin{align*}
  \int_\R \left(\partial_x \udx(t,x)\right)^2 \dx
  & = \Dx\sum_j\left(q_j(t)\right)^2 \\ & 
  \le \frac{1}{\Dt}\left(\left(t-t^n\right)
    \Dx\sum_j \left(q^{n+1}_j\right)^2 
    +\left(t^{n+1}-t\right)\Dx
    \sum_j \left(\qnj\right)^2 \right).
\end{align*}
Hence, 
$$
\norm{\udx(t,\cdot)}_{H^1(\R)}^2 \le
\frac{1}{\Dt}\left(
\left(t-t^n\right)\norm{u^{n+1}}_{H^1(\R)}^2 +
\left(t^{n+1}-t\right)\norm{u^{n}}_{H^1(\R)}^2\right).
$$
Let $s\in I^m$ with $m\le n$. Then,  using \eqref{eq:enest},
\begin{align*}
  \norm{\udx(t,\cdot)}_{H^1(\R)}^2
  &\le 
  \frac{1}{\Dt}\left(
    \left(t-t^n\right)\norm{u^{n+1}}_{H^1(\R)}^2 +
    \left(t^{n+1}-t\right)\norm{u^{n}}_{H^1(\R)}^2\right)\\
  &\le e^{\Dt(n-m)\Dx^\theta} \norm{\udx(s,\cdot)}_{H^1(\R)}^2 \\
  &\quad + 
  \frac{e^{\Dt(n-m)\Dx^\theta}}{\Dt}\biggl[
    \left(t-t^n\right)
    \left(\norm{u^{m+1}}_{H^1(\R)}^2-\norm{\udx(s,\cdot)}_{H^1(\R)}^2\right) \\
    &\hphantom{\quad+\frac{1}{\Dt}\biggl[}\quad + 
    \left(t^{n+1}-t\right)
    \left(\norm{u^{m}}_{H^1(\R)}^2-\norm{\udx(s,\cdot)}_{H^1(\R)}^2\right)
  \biggr]\\
  &\le  e^{\Dt(n-m)\Dx^\theta} \norm{\udx(s,\cdot)}_{H^1(\R)}^2 
  \\&\qquad 
  + e^{\Dt(n-m)\Dx^\theta}
  \abs{\norm{u^{m+1}}_{H^1(\R)}^2-\norm{u^{m}}_{H^1(\R)}^2}\\
  &\le e^{\Dt(n-m)\Dx^\theta} \norm{\udx(s,\cdot)}_{H^1(\R)}^2 \\
  &\qquad  + 
  e^{\Dt(n-m)\Dx^\theta}
  \left(e^{\Dt\Dx^\theta}-1\right)\norm{u^{m}}_{H^1(\R)}^2.  
\end{align*}
This implies \eqref{eq:uweakstar} and \eqref{eq:hdecrease}.
  
Next we prove that $\seq{\udx}_{\Dx>0}$ is uniformly bounded 
in $W^{1,2+\alpha}((0,T)\times (a,b))$. We can assume that 
$a-1=x_{j_a}$ and $b+1=x_{j_b}$ for some integers $j_a$ and $j_b$. 
Since $q\mapsto \abs{q}^{2+\alpha}$ is convex,
\begin{equation}\label{eq:udxhigh}
        \begin{aligned}
                & \int_0^T\int_a^b \abs{\partial_x\udx}^{2+\alpha}\dx \\ 
                & \quad \le \sum_n\Dx\sum_j\left( \left(t^{n+1}-t\right)\abs{\qnj}^{2+\alpha}
                +\left(t-t^n\right)\abs{\qnj}^{2+\alpha}\right)\le C,
        \end{aligned}
\end{equation}
for some constant $C=C(\alpha,a,b,u_0)$, where we have 
also used \eqref{eq:highint}.

Now set $\sigma=(x-x_{j-1/2})/\Dx$. Then, for $x\in I_j$, we have
\begin{align*}
        \abs{\partial_t \udx(t,x)}
        & =\abs{(1-\sigma)\DT\unjmh+\sigma\DT\unjph}
        \\ & \le (1-\sigma)\abs{\DT\unjmh}+\sigma\abs{\DT\unjph}.
\end{align*}
Furthermore, by the uniform bounds on $\Dp\Pnj$ and $\udx$,
$$
\abs{\DT\unjph}\le C(1+\abs{\qnj}+\abs{q^n_{j-1}}).
$$
Using this,
$$
\int_0^T\int_a^b \abs{\partial_t \udx}^{2+\alpha}\dx\le 
C\left(1+\Dt\sum_{n=0}^N\Dx\sum_{{j_a}}^{{j_b}} 
\abs{\qnj}^{2+\alpha}\right)\le C.
$$
Now $\seq{\udx}_{\Dx>0}\subset W^{1,2+\alpha} \subset \subset C^{0,\el{}}$ 
on $(0,T)\times (a,b)$ with $\ell=1-2/(2+\alpha)$. 
Therefore, along a subsequence, $\udx\to u$ uniformly 
in $(0,T)\times (a,b)$ as $\Dx\to 0$.

Let us show that the limit satisfies the initial condition \eqref{eq:uinitialconv}. 
Fix $\bar{x}\in\R$ and let $t\in(0,1)$. We have $\bar{x}\in I_j$ for 
some $j$ and $\udx(x_{j-1/2},0)=u_0(x_{j-1/2})$, so that
\begin{align*}
        \abs{\udx(0,\bar{x})-u_0(\bar{x})}
        &\le \abs{\udx(0,\bar{x})-\udx(0,x_{j-1/2})}
        +\abs{u_0(x_{j-1/2})-u_0(\bar{x})}\\
        & \le C\left(\bar{x}-x_{j-1/2}\right)^\ell.
\end{align*}
Consequently,
\begin{align*}
        \abs{u(t,\bar{x})-u_0(\bar{x})}
        & \le \abs{u(t,\bar{x})-\udx(t,\bar{x})}\\ 
        &\qquad +\abs{\udx(t,\bar{x})-\udx(0,\bar{x})}
        +\abs{\udx(0,\bar{x})-u_0(\bar{x})}\\
        &\le \abs{u(t,\bar{x})-\udx(t,\bar{x})} + Ct^\ell + C\Dx^\ell.
\end{align*}
Now we can let $\Dx\to 0$ and then $t\to 0$ to conclude that $u(\bar{x},0)=u_0(\bar{x})$. 
This draws to a close the proof of the lemma.
\end{proof}

\begin{lemma}\label{lem:Pconvergence} 
There exists a limit function
$$
P\in L^\infty(0,T;W^{1,\infty}(\R))\cap L^{\infty}(0,T;W^{1,1}(\R)),
$$
such that along a subsequence as $\Dx\to 0$
\begin{equation}\label{eq:Pconvergence}
        \text{$\Pdx\to P$ in $L^p_{\loc}(\Do)$}, \quad 1\le p<\infty.
\end{equation}
\end{lemma}

\begin{proof}
By the bounds on $\Pnj$ in Lemma~\ref{lem:boundP}, we see that $\Pdx$ 
is bounded in $L^\infty$ uniformly in $\Dx$. Next we show 
that $\seq{\partial_t \Pdx}_{\Dx>0}$ is bounded in 
$L^1(\Do)$. For $t\in [t^n,t^{n+1})$ and $x\in I_{j+1/2}$,
$$
\partial_t\Pdx(t,x)= \DT\Pnj + (x-x_j)\DT\Dp \Pnj =
\left(1-\sigma\right)\DT\Pnj + \sigma \DT P^{n}_{j+1},
$$
where $\sigma=(x-x_j)/\Dx$. Write $\DT\Pnj=X^n_j+Y^n_j$, 
where $X^n_j$ and $Y^n_j$ solve
\begin{align*}
        X^n_j - \Dm\Dp X^n_j &= 
        \DT\left( \left(\unjmh\wedge 0\right)^2
        +\left(\unjph\vee 0\right)^2 \right),
        \\
        Y^n_j -\Dm\Dp Y^n_j &=
        \frac{1}{2} \DT\left(\left(\qnj\right)^2\right).
\end{align*}
  
Then $\norm{X^n}_{\el1}$ is bounded by the $\el1$ norm of the 
corresponding right hand side avbove. By the discrete chain rule
$$
\DT\left(\unjmh\wedge 0\right)^2 \le 2\left(\unjmh\wedge 
0\right)\DT\unjmh + \Dt \left(\DT\unjmh\right)^2.
$$
Estimating the first term here
\begin{align*}
        & \left(\unjmh \wedge 0\right)\DT\unjmh
        \\ &\qquad
        = -\left(\unjmh\wedge 0\right)\left[\left(\unjmh\vee 0\right)q^n_{j-1}
        + \left(\unjmh\wedge 0\right)\qnj +\Dp\Pnj\right]
        \\ & \qquad
        =-\left(\unjmh\wedge 0\right)
        \left[\left(\unjmh\wedge 0\right)\qnj+\Dp\Pnj\right].
\end{align*}
This means that
$$
\norm{\left(u^n\wedge 0\right)\DT u^n}_{\el1}
\le C\Biggl[ \norm{u^n}_{\el2}\norm{q^n}_{\el2}
+\norm{\Dp P^n}_{\el1}\Biggr]\le C.
$$
Similarly, we have
\begin{align*}
        \left(\DT\unjph\right)^2& \!\le\! C\left[\left(\unjph\vee 0\right)^2
        +\left(\qnj\right)^2+\left(\unjph\wedge 0\right)^2
        +\left(q^n_{j+1}\right)^2 +\left(\Dp\Pnj\right)^2\right]
        \\ & \le C \left[ \left(\unjph\right)^2+\left(\qnj\right)^2+\left(q^n_{j+1}\right)^2
        +\abs{\Dp\Pnj}\right],
\end{align*}
and thus $\norm{\left(\DT u^n\right)^2}_{\el1}\le C$. We have shown that
\begin{equation}\label{eq:firstDtbnd}
        \norm{X^n}_{\el1}\le \norm{\seq{\DT\left[(u^n_{j-1/2}\wedge 0)^2
        +(u^n_{j+1/2}\vee 0)^2\right]}_j}_{\el1} \le C.
\end{equation}
Next, using \eqref{eq:fq-disc-II} with $f(q)=q^2/2$ we have that
\begin{equation}\label{eq:qsquarescheme}
        \begin{aligned}
                \frac{1}{2}\DT\left(\qnj\right)^2 &
                =\underbrace{-\Dm\left[\left(\unjph\vee 0\right)
                \left(\qnj\right)^2\right]+\Dp\left[\left(\unjmh\wedge 0\right)
                \left(\qnj\right)^2\right]}_{a^n_j}
                \\ &\qquad\qquad\quad
                + \underbrace{A^n_j\qnj
                -I_{\Dx,f'',j}+\Dt\left(\DT\qnj\right)^2}_{b^n_j},
        \end{aligned}
\end{equation}
where $A^n_j$ is defined in \eqref{eq:Anj}  and
$$
I_{\Dx,f'',j}= \Dx\left\{\left(\unjmh\vee 0\right)\left(\Dm\qnj\right)^2
-\left(\unjph\wedge 0\right)\left(\Dp\qnj\right)^2\right\}\ge 0.
$$
We write $Y^n_j=Y^{a,n}_j+Y^{b,n}_j$ where
$$
Y^{a,n}_j=\left(I-\Dm\Dp\right)^{-1} a^n_j\quad\text{and}\quad
Y^{b,n}_j=\left(I-\Dm\Dp\right)^{-1} b^n_j.
$$
Now $\norm{Y^{b,n}}_{\el1}\le \norm{b^n}_{\el1}$, and therefore we compute
\begin{align*}
        \norm{A^n q^n}_{\el1}&
        =\Dx\sum_j \abs{A^n_j\qnj}
        \le C\norm{q^n}_{\el2}
        \left( \Dx\sum_j \left[ \left(\Pnj\right)^2
        +\left(\unjph\right)^2\right]\right)^{1/2}\\ 
        & \le C\norm{q^n}_{\el2}\left(\Dx\sum_j \left[\abs{\Pnj}
        +\left(\unjph\right)^2\right]\right)^{1/2} \le C.
\end{align*}
By \eqref{eq:dtsbound}, the $L^1$ norm of $\Dt(\DT\qnj)^2$ is of 
the same order as $\Dx^\theta$. Then, summing \eqref{eq:qsquarescheme} 
over $n$ and $j$, we arrive at
\begin{equation*}
  \Dt\Dx\sum_{n,j} I_{\Dx,f'',j}
  \le CT+\Dx\sum_j\left[ \left(q^N_j\right)^2 + \left(q^0_j\right)^2\right]
  +\mathcal{O}\left(\Dx^\theta\right)\le C.
\end{equation*}
This means that
\begin{equation}\label{eq:Ybbnd}
        \Dt\Dx \sum_{n,j}\abs{Y^{b,n}_j}\le C.
\end{equation}
  
Now let
$$
L_j = h \sum_i e^{-\kappa\abs{i-j}} \Dpm K_i=
-h\sum_i\Dmp\bigl( e^{-\kappa\abs{i-j}}\bigr) K_i,
$$
for some sequence $\seq{K_j}_j\in \el1$. Since
$$
\abs{\Dpm e^{-\kappa\abs{i-j}}}\le Ce^{-\kappa\abs{i-j}},
$$
we get
$$
\norm{L}_{\el1}=\Dx\sum_j \abs{L_j} \le 
\Dx h C \sum_{i,j}e^{-\kappa\abs{i-j}} 
\abs{K_i} \le C\norm{K}_{\el1}.
$$
Using this,
$$
\norm{Y^{a,n}}_{\el1}\le C\norm{q^n}_{\el2}^2.
$$
Combining this with \eqref{eq:Ybbnd} and \eqref{eq:firstDtbnd} 
we see that
$$
\Dt\Dx\sum_{n,j} \abs{\DT\Pnj} \le C,
$$
and therefore
$$
\int_0^T \int_{\R} \abs{\partial_t \Pdx}\dx\dt \le C.
$$
Hence $\seq{\Pdx}_{\Dx>0}$ is bounded in $W^{1,1}(\Do)$. 
Combining this with the $L^\infty$ estimates found 
in Lemma~\ref{lem:boundP} yields the existence of a 
convergent subsequence as claimed in \eqref{eq:Pconvergence}.
\end{proof}

\section{Strong convergence result}
\label{sec:strongconv}

We now show that the sequence $\seq{\qdx}_{\Dx>0}$, cf.~\eqref{eq:qdxdef}, has 
a strongly converging subsequence. 
This result is a key point of the convergence analysis.

\begin{lemma}\label{lem:boundq2}
Fix $1\le p <3$ and $1\le r <1+\frac{\theta}{2}$. Then there exist 
two functions $q\in L^{p}_{\loc}(\Do)$, $\wlim{q^2}\in L^{r}_{\loc}(\Do)$ 
such that for a subsequence as $\Dx\to 0$,
\begin{align}
  &\text{$\qdx\weakstar q$ in $L^{\infty}(0,T;L^{2}(\R))$},
  &\text{$\qdx \weak q$ in $L^{p}_{\loc}(\Do)$}, \label{eq:q1*}\\ 
  & \text{$\qdx^2 \weak \wlim{q^2}$ in $L^{r}_{\loc}(\Do)$}, \label{eq:q2}
\end{align}
for all $a,b\in \R$, $a<b$. Moreover,
\begin{equation}\label{eq:q3} 
  \text{$q^2(t,x)\le \wlim{q^2}(t,x)$ for a.e.~$(t,x)\in \Do$}
\end{equation}
and
\begin{equation}\label{eq:der1}
  \px u=q \quad \text{in the sense of distributions on $\Do$.}
\end{equation}
Finally, there is a positive constant $C$ such that
\begin{equation}\label{eq:oleinik_limit}
  q(t,x)\le  \frac{2}{t}+C\norm{u_0}_{H^1(\R)}, \qquad 
  t\in (0,T),\, x\in \R.
\end{equation}
\end{lemma}

\begin{proof} 
Claims \eqref{eq:q1*}, \eqref{eq:q2} are direct consequences of 
Lemmas \ref{lem:enest} and \ref{lem:highint}.  
Claim \eqref{eq:q3} is true thanks to \eqref{eq:q2} and the 
convexity of $g=q^2$, cf.~Lemma \ref{lem:prelim}, while 
\eqref{eq:der1} is a consequence of the definitions of $\qdx$ 
and $\udx$, cf.~ \eqref{eq:qdxdef} and \eqref{eq:udxdef}.
  
We conclude by proving \eqref{eq:oleinik_limit}. Fix $t>0$, and let
$\Dt$ be so small that $t\in I^n$ with $n>0$. Not that $n\to \infty$
as $\Dt\to 0$, and $t^n$ and $t^{n+1}$ both tend to $t$.
From the definition of $\qdx$ and \eqref{eq:oleinikbnd} we have that
\begin{equation}\label{eq:olfin1}
  \begin{split}
    \qdx(t,x)&=\qnj+(t-t^n)\DT \qnj \\ 
    &=\frac{t^{n+1}-t}{\Dt}\qnj+\frac{t-t^n}{\Dt}q^{n+1}_j \\ 
    &\le\frac{t^{n+1}-t}{\Dt} \left(\frac{2}{t^n}+\hat{C}\right)
    +\frac{t-t^n}{\Dt}\left(\frac{2}{t^n}+\hat{C}\right)
    \\ & =\frac{2}{t}+\hat{C}+2f_{\Dt}(t),
  \end{split}
\end{equation}
where $\hat{C}:=C\norm{u_0}_{H^1(\R)}$ and
for every  $t\in [t^n,t^{n+1})$, $x\in I_j$, with
$$
f_{\Dt}(t) =\frac{t-t^n}{\Dt}\frac{1}{t^{n+1}}
+\frac{t^{n+1}-t}{\Dt}\frac{1}{t^{n}}-\frac{1}{t}.
$$
Observe that
$$
f_{\Dt}'(t)=\frac{1}{\Dt}\left(\frac{1}{t^{n+1}}-\frac{1}{t^n}\right)
+\frac{1}{t^2}=-\frac{1}{t^nt^{n+1}}+\frac{1}{t^2},
$$
so $f_{\Dt}'(t)=0\Longleftrightarrow t=\sqrt{t^nt^{n+1}}\in (t_n,t_{n+1})$, and 
in particular 
\begin{align*}
  \sup\limits_{t\in [t^n,t^{n+1})} f_{\Dt}(t)
  &= f_{\Dt}(\sqrt{t^nt^{n+1}})
  =\frac{\left(\sqrt{t^{n+1}}-\sqrt{t^n}\right)^2}{t^nt^{n+1}}
  \\ & \le \left(\frac{\Dt}{2\sqrt{t^n}}\right)^2
  \frac{1}{t^n t^{n+1}} \le \frac{\Dt^2}{8(t^n)^2t^{n+1}}\to 0.
\end{align*}
Therefore \eqref{eq:oleinik_limit} follows 
from \eqref{eq:q1*} and \eqref{eq:olfin1}.
\end{proof}

In view of the weak convergences stated in \eqref{eq:q1*}, we have
that for any function $f\in C^1(\R)$ with $f'$ bounded
\begin{equation}\label{eq:weak-f-limits}
        \begin{split} 
                & f(\qdx) \weakstar \wlim{f(q)} 
                \quad \text{in $L^{\infty}(0,T;L^{2}(\R))$}, 
                \\ & \text{$f(\qdx) \weak 
                \wlim{f(q)}$ in $L^{p}_{\loc}(\Do)$, $1\le p<3$},
  \end{split}
\end{equation} 
where the same subsequence of $\Dx\to 0$ applies to 
any $f$ from the specified class.

In what follows, we let $\wlim{qf(q)}$ and $\wlim{f'(q)q^2}$ denote
the weak limits of $\qdx f(\qdx)$ and $f'(\qdx)\qdx^2$, respectively,
in $L^{r}_{\loc}(\Do)$, $1\le r<\frac{3}{2}$.

\begin{lemma} \label{lem:weak-f(q)}
For any convex function $f\in C^1(\R)$ with $f'$ 
bounded we have that
\begin{equation*}
        \begin{split}
                \iint_{\Do} & \left(\wlim{f(q)} \pt \test + u \wlim{f(q)}
                \px \test\right) \dx\dt + \int_{\R} f(q_0(x)) \test(0,x)\dx
                \\ &\qquad \ge \iint_{\Do} \left(\frac{1}{2}\wlim{f'(q)q^2}
                -\wlim{qf(q)} + \left(P-u^2\right) \wlim{f'(q)}\right)\test\dx\dt,
        \end{split}
\end{equation*}
for any nonnegative $\test \in C^\infty_c(\hDo)$.
\end{lemma}

\begin{proof} 
Set
\begin{equation}\label{eq:testdef}
        \test_j(t)=\frac{1}{\Dx}\int_{I_j}\test(t,x)\dx,\quad
        \test_j^n=\frac{1}{\Dx\Dt}\iint_{I_j^n}\test(t,x)\dx\dt.
\end{equation} 
We multiply \eqref{eq:fq-disc-II} by $\Dx\Dt\test_j^n$, sum over $n,j$, and 
take into account the convexity of $f$. After partial summations, the final result reads
\begin{equation}\label{eq:wf0}
        E_0+E_1+E_2+E_3+E_4+E_5\ge 0,
\end{equation}
where
\begin{align*}
        E_0=& \Dx \sum_{j}  f(q_j^0) \test_j^0,\\ 
        E_1=& \Dt\Dx\sum\limits_{n,j}f(\qnj)D_t^-\test_j^n,\\
        E_2=& \Dt\Dx\sum\limits_{n,j}\left[\left(\unjph\vee 0\right)
        f(\qnj)\Dp\test_j^n+\left(\unjmh\wedge 0\right)f(\qnj)\Dm\test_j^n\right],\\ 
        E_3=&-\Dt\Dx\sum\limits_{n,j}\left[\frac{(\qnj)^2}{2}f'(\qnj)-\qnj f(\qnj)\right]\test_j^n,\\
        E_4=&-\Dt\Dx\sum\limits_{n,j}A_j^n f'(\qnj)\test_j^n,\\
        E_5=&\Dt^2\Dx\sum\limits_{n,j} f''\left(q^{n+1/2}_j\right)
        \left(\DT\qnj\right)^2\test_j^n.
\end{align*}

By \eqref{eq:init-q}, 
\begin{equation*}
        \norm{q_0-q_{\Dx,0}}_{L^2(\R)}
        \to 0 \quad \text{as $\Dx\to 0$,} 
\end{equation*}
where $q_0:=\px u_0$ and $q_{\Dx,0}:= \px \udx|_{t=0}$, so that 
$$
\Dx \sum_{j}  f(q_j^0) \test_j^0\to  \int_{\R} f(q_0(x)) \test(0,x)\dx 
\quad \text{as $\Dx\to 0$}.
$$

We split $E_1$ in three parts:
\begin{equation}\label{eq:wf1}
        E_1=E_{1,1}+E_{1,2}+E_{1,3},
\end{equation}
where
\begin{align*}
        E_{1,1}= &\iint_{\Do}f(\qdx)\pt\test\dx\dt,\\
        E_{1,2}= &\sum\limits_{n,j}\iint_{I_j^n}\left(f(\qnj)-f(\qdx)\right)D_t^-\test_j^n \dx\dt\,\\
        E_{1,3}= &\sum\limits_{n,j}\iint_{I_j^n}f(\qdx)\left(D_t^-\test_j^n-\pt\test\right)\dt.
\end{align*}
Due to \eqref{eq:weak-f-limits}
\begin{equation}\label{eq:wf2}
        E_{1,1}\to \iint_{\Do} \wlim{f(q)} \pt \test \dx\dt.
\end{equation}
Due to the boundedness of $f'$,
\begin{equation}\label{eq:boundf}
        \abs{f(\qnj)-f(\qdx)}\le c_1(t-t^n)|\DT\qnj|,
\end{equation}
for each $(t,x)\in [t^n,t^{n+1})\times I_j$, where $c_1>0$ is a finite constant. 
In view of \eqref{eq:dtsbound},
\begin{equation}\label{eq:wf3}
        \begin{split}
                \abs{E_{1,2}}&\le c_1\Dx\sum\limits_{n,j}
                \abs{D_t^-\test_j^n}|\DT\qnj|\int_{t^n}^{t^{n+1}}(t-t^n)\dt \\&
                \le c_1\Dx\Dt^2\sum\limits_{n,j}\abs{D_t^-\test_j^n}|\DT\qnj|\le
                c_1\Dt\norm{D_t^-\test}_{\el2} \norm{\DT q}_{\el2} \to 0,
        \end{split}
\end{equation}
as $\Dx\to 0$. Finally, since
$$
\abs{D_t^-\test-\pt\test}\le c_2\Dx,
$$
for some constant $c_2>0$, we have that
\begin{equation}\label{eq:wf4}
        \abs{E_{1,3}}\le c_2\Dx\iint_{\mathrm{supp}(\test)}\abs{f(\qdx)}\dx\dt\to 0
        \quad \text{as $\Dx\to 0$.}
\end{equation}
Clearly \eqref{eq:wf1}, \eqref{eq:wf2}, \eqref{eq:wf3}, and \eqref{eq:wf4} imply
\begin{equation}\label{eq:wf4'}
        E_{1}\to \iint_{\Do} \wlim{f(q)} \pt \test \dx\dt \quad 
        \text{as $\Dx\to 0$.}
\end{equation}

Next, we split $E_2$ into four parts:
\begin{equation}\label{eq:wf1.1}
        E_2=E_{2,1}+E_{2,2}+E_{2,3}+E_{2.4},
\end{equation}
where
\begin{align*}
        E_{2,1} = &\iint_{\Do}\udx f(\qdx)\px\test\dx\dt, \\
        E_{2,2} = &\sum\limits_{n,j}\iint_{I_j^n}
        \Biggl [ \left( (\unjph\vee 0)- (\udx\vee 0)\right) f(\qnj)\Dp\test_j^n 
        \\ & \qquad \qquad\qquad
        + \left((\unjmh\wedge 0)-(\udx\wedge0) \right)f(\qnj)\Dm\test_j^n\Biggr]\dx \dt,\\
        E_{2,3}=&\sum\limits_{n,j}\iint_{I_j^n}
        \Biggl[ (\udx\vee 0)\left(f(\qnj)-f(\qdx)\right)\Dp\test_j^n
        \\& \qquad\qquad\qquad
        +(\udx\wedge 0) \left(f(\qnj)-f(\qdx)\right)\Dm\test_j^n\Biggr]\dx\dt,\\
        E_{2,4}=&\sum\limits_{n,j}
        \iint_{I_j^n}\Biggl[\left(\udx\vee0\right)f(\qdx)
        \left(\Dp\test_j^n-\px\test\right)\\&
        \qquad\qquad\qquad
        +(\udx\wedge 0)f(\qdx) \left(\Dm\test_j^n-\px\test\right)\Biggr]\dx\dt.
\end{align*}
Due to \eqref{eq:weak-f-limits},
\begin{equation}\label{eq:wf2.1}
        E_{2,1}\to \iint_{\Do} \wlim{qf(q)} \px \test \dx\dt.
\end{equation}
Using the definition of $\udx$,
\begin{equation*}
        \begin{split}
                \abs{\unjmh-\udx} & \le 
                \Dx\left(\abs{\qnj}+\abs{q_j^{n+1}}\right),\\
                \abs{\unjmh-\udx} & \le
                \Dx\qnj\left(\abs{\qnj}+\abs{q_j^{n+1}}\right)
                +\abs{\unjph-\unjmh},
        \end{split}
\end{equation*}
so, by Lemmas \ref{lem:enest} and \ref{lem:convlemma1},
\begin{equation}\label{eq:wf3.1.1}
        \begin{split}
                \abs{E_{2,2}} &\le \Dx^2\Dt\sum\limits_{n,j}
                \left(\abs{\qnj}+\abs{q_j^{n+1}}\right)
                \abs{f(\qnj)}\left(\abs{\Dp\test_j^n}+\abs{\Dm\test_j^n}\right)
                \\ & \qquad\qquad\quad
                + \Dx\Dt\sum\limits_{n,j}\abs{\unjph-\unjmh}
                \abs{f(\qnj)}\abs{\Dp\test_j^n} \\ & 
                \le 2\Dx\norm{\Dp\test}_{\el\infty}\norm{q}_{\el2}
                \norm{f(q)}_{\el 2}
                \\ &\qquad\qquad\quad
                + \norm{f(q)}_{\el2}\norm{\left\{ \left(\unjph-\unjmh\right)
                \Dp\test_j^n\right\}_{n,j}}_{\el2}\to 0,
        \end{split}
\end{equation}
as $\Dx\to 0$. Using \eqref{eq:boundf} and \eqref{eq:dtsbound}, we deduce
\begin{equation}\label{eq:wf3.1}
  \begin{split}
    \abs{E_{2,3}} & \le c_1\Dx\sum\limits_{n,j}
    \left(\abs{\Dm\test_j^n}+\abs{\Dp\test_j^n}\right)
    \\ & \qquad\qquad\quad \times
    \left(\abs{\unjph}+\abs{\unjmh}\right)
    \int_{t^n}^{t^{n+1}}(t-t^n) \dt\\
    & = c_1\Dx\Dt^2\sum\limits_{n,j}
    \left(\abs{\Dm\test_j^n}+\abs{\Dp\test_j^n}\right)
    \left(\abs{\unjph}+\abs{\unjmh}\right)
    \\ & \le 4c_1\Dt\norm{\Dp\test}_{\el 2} \norm{u}_{\el 2}\to 0,
  \end{split}
\end{equation}
as $\Dx\to 0$. Finally, since $\abs{D_\pm\test-\pt\test}\le c_3\Dx$, for 
some constant $c_3>0$, we obtain
\begin{equation}
  \label{eq:wf4.1}
  \abs{E_{2,4}}\le c_3\Dx
  \iint_{\mathrm{supp}(\test)}
  \abs{\qdx}\abs{f(\qdx)}\dt\dx\to 0
  \quad \text{as $\Dx\to 0$.}
\end{equation}
Note that \eqref{eq:wf1.1}, \eqref{eq:wf2.1}, \eqref{eq:wf3.1.1},
\eqref{eq:wf3.1}, and \eqref{eq:wf4.1} imply
\begin{equation}\label{eq:wf4''}
        E_{2}\to \iint_{\Do} \wlim{qf(q)} \px \test \dx\dt  \quad \text{as $\Dx\to 0$.}
\end{equation}

We split $E_3$ into three parts:
\begin{equation}\label{eq:wf1.2}
        E_3=E_{3,1}+E_{3,2}+E_{3,3},
\end{equation}
where
\begin{align*}
        E_{3,1} & =-\iint_{\Do}
        \left[\frac{(\qdx)^2}{2}f'(\qdx)-\qdx f(\qdx)\right]\test\dx\dt,
        \\ E_{3,2}& = -\sum\limits_{n,j}\iint_{I_j^n}\Biggl[\frac{(\qnj)^2}{2}f'(\qnj)-\qnj f(\qnj)
        \\ & \qquad\qquad \qquad\qquad \quad
        - \frac{(\qdx)^2}{2}f'(\qdx)+\qdx f(\qdx)\Biggr]\test_j^n\dx\dt, \\
        E_{3,3} &= -\sum\limits_{n,j}\iint_{I_j^n}\left[\frac{(\qdx)^2}{2}f'(\qdx)-\qdx f(\qdx)\right]
        \left(\test_j^n-\test\right)\dx\dt,
\end{align*}
Due to \eqref{eq:weak-f-limits},
\begin{equation}\label{eq:wf2.7}
        E_{3,1}\to -\iint_{\Do}
        \left(\frac{\wlim{q^2f'(q)}}{2}-\wlim{qf(q)}\right) \test \dx\dt
        \quad \text{as $\Dx\to 0$.}
\end{equation}
Using the boundedness of $f'$, we can estimate as follows:
\begin{align*}
        &\abs{\frac{(\qnj)^2}{2}f'(\qnj)
        -\qnj f(\qnj)-\frac{(\qdx)^2}{2}f'(\qdx)+\qdx f(\qdx)} 
        \\& \le \frac{\abs{(\qnj)^2-\qdx^2}}{2} |f'(\qnj)|
        +\frac{(\qdx)^2}{2}\abs{f'(\qnj)-f'(\qdx)} 
        \\ & \le \frac{\abs{\qnj-\qdx}(\abs{\qnj}
        +\abs{\qdx})}{2}\norm{f'}_{L^\infty}
        +\frac{(\qdx)^2}{2}\norm{f''}_{L^\infty}\abs{\qnj-\qdx} 
        \\& \le \frac{(t-t^n)\abs{\DT \qnj}(\abs{\qnj}
        +\abs{\qdx})}{2}\norm{f'}_{L^\infty}
        +\frac{(\qdx)^2}{2}\norm{f''}_{L^\infty}
        (t-t^n)\abs{\DT \qnj}.
\end{align*}
Hence, taking into account \eqref{eq:dtsbound} and \eqref{eq:highint},
\begin{equation}\label{eq:wf2.8}
  \begin{split}
    \abs{E_{3,2}}&\le \Dx\sum\limits_{n,j}\int_{I^n}
    \Biggl[\frac{(t-t^n)\abs{\DT \qnj}
      (\abs{\qnj}+\abs{\qdx})}{2}\norm{f'}_{L^\infty}
    \\ & \qquad\qquad\qquad\qquad\qquad 
    +\frac{(\qdx)^2}{2}\norm{f''}_{L^\infty}(t-t^n)\abs{\DT
      \qnj}\Biggr]\test_j^n\dt 
    \\ &\le \Dx\Dt^2\sum\limits_{n,j}\left[\frac{\abs{\DT \qnj}(\abs{\qnj}
        +\abs{\qdx})}{2}\norm{f'}_{L^\infty}
      +\frac{(\qdx)^2}{2}\norm{f''}_{L^\infty}\abs{\DT \qnj}\right]\test_j^n
    \\ &\le \Dt\Biggl[\frac{\norm{\DT q}_{\el 2}(\norm{q^n}_{\el 2}
      +\norm{\qdx}_{L^2(\R)})}{2}\norm{f'}_{L^\infty}
    \\ & \qquad\qquad\qquad\qquad\qquad
    +\mathcal{O}\left(\frac{1}{\Dx}\right)\norm{u}_{\el\infty}
    \norm{q}_{\el2}\norm{\DT q^n}_{\el2} \Biggl]
    \norm{\test}_{L^\infty}
    \\ &\le \Dt\left[\mathcal{O}\left(\sqrt{\frac{\Dx^\theta}{\Dt}}
      +\sqrt{\frac{\Dx^{\theta}}{\Dt}}\frac{1}{\Dx}\right)
    \right]=
    \mathcal{O}\left(\Dx^\theta\right)\to 0 \quad \text{as $\Dx\to 0$.}
  \end{split}
\end{equation}
Since $\abs{\test_j^n-\test}={\mathcal O}(\Dx)$,
\begin{equation}\label{eq:wf2.8.0}
        \abs{E_{3,3}}\le{\mathcal O}(\Dx) \iint_{\supp (\test)}
        \abs{\frac{(\qdx)^2}{2}f'(\qdx)-\qdx f(\qdx)}\dx\dt\to 0
        \quad \text{as $\Dx\to 0$}.
\end{equation}
We have that \eqref{eq:wf1.2}, \eqref{eq:wf2.7}, \eqref{eq:wf2.8}, and
\eqref{eq:wf2.8.0} imply
\begin{equation}\label{eq:wf4,0,0}
        E_{3}\to \iint_{\Do} \wlim{qf(q)} \px \test \dx\dt  \quad \text{as $\Dx\to 0$.}
\end{equation}

We split the term $E_4$ into three parts:
\begin{equation*}
        E_4=E_{4,1}+E_{4,2}+E_{4,3},
\end{equation*}
where
\begin{align*}
        E_{4,1}&=-\iint_{\Do}A_{\Dx}f'(\qdx)\test\dx\dt,
        \\ E_{4,2}&=-\sum\limits_{n,j}\iint_{I_j^n}
        \left(A_j^nf'(\qnj)-A_{\Dx}f'(\qdx)\right)\test_j^n\dx\dt,
        \\ E_{4,3}&=-\sum\limits_{n,j}\iint_{I_j^n}
        A_{\Dx}f'(\qdx)\left(\test_j^n-\test\right)\dx\dt,
\end{align*}
where $A_{\Dx} = P_{\Dx} - \left(\udx\right)^2$. 
Lemmas \ref{lem:convlemma1} and 
\ref{lem:Pconvergence}, cf.~also \eqref{eq:weak-f-limits}, imply that
\begin{equation*}
        E_{4,1}\to - \iint_{\Do} \left(P-u^2\right) \wlim{f'(q)}\test \dx\dt
        \quad \text{as $\Dx\to 0$}.
\end{equation*}
Continuing, it is not hard to see that 
$\abs{E_{4,2}}=\mathcal{O}(\Dx)\to 0$ as $\Dx\to 0$.
Moreover, since $\test_j^n-\test=\mathcal{O}(\Dx)$, 
$E_{4,3}\to 0$ as $\Dx\to 0$. 
Summarizing,
\begin{equation}\label{eq:wf4.5}
        E_{4}\to - \iint_{\Do} \left(P-u^2\right)\wlim{f'(q)}\test \dx\dt
        \quad \text{as $\Dx\to 0$.}
\end{equation}

Finally, regarding $E_5$, due to \eqref{eq:dtsbound}, we 
conclude that as $\Dx\to 0$
\begin{equation}\label{eq:wf5}
        \abs{E_{5}}\le \norm{\test}_{L^\infty(\Do)}\Dt^2\Dx
        \sum\limits_{n,j}f''\left(q^{n+1/2}_j\right)
        \left(\DT\qnj\right)^2\test_j^n\to 0.
\end{equation}

The lemma now follows from \eqref{eq:wf0}, \eqref{eq:wf4'}, \eqref{eq:wf4''}, 
\eqref{eq:wf4,0,0}, \eqref{eq:wf4.5}, and \eqref{eq:wf5}.
\end{proof}

We know that $\seq{(\qdx)^2}_{\Dx>0}$ is bounded in
$L^\infty(0,T;L^1(\R)) \cap L_{\loc}^r(\Do)$, for any $1\le r<1+\theta/2$.
Additionally, using \eqref{eq:fq-disc-II} with $f(q)=\tfrac{q^2}{2}$,
we can show that the mapping $t\mapsto\int_{\R} (\qdx)^2 \test \dx$ is
equi-continuous on $[0,T]$, for every $\test \in C^\infty_c(\R)$.
Hence, in view of Lemma \ref{lem:timecompactness},
\begin{equation}\label{eq:weak_uniform_conv-q}
        \int_{\R}  (\qdx)^2 \test \dx \to \int_{\R} 
        \wlim{q^2} \test \dx  
        \quad \text{uniformly on $[0,T]$,}
\end{equation}
and 
\begin{equation}\label{eq:weakcont-q}
        t \mapsto \int_{\R}  \wlim{q^2}\test \dx
        \quad \text{is continuous on $[0,T]$.}
\end{equation}
The statements \eqref {eq:weak_uniform_conv-q} and \eqref{eq:weakcont-q} 
hold with $\qdx^2$, $\wlim{q^2}$ replaced
respectively by $f(\qdx)$, $\wlim{f(\qdx)}$, for any convex
function $f\in C^1(\R)$ with $f'$ bounded.

\begin{lemma}\label{lem:weak-q}
Let $q$ and $\wlim{q^2}$ be the weak limits 
identified in Lemma \ref{lem:boundq2}. Then
\begin{equation}\label{eq:weak-q}
        \begin{split}
                & \iint_{\Do} \left( q \pt \test + u q \px \test\right) \dx\dt
                +\int_{\R} q_0(x) \test(0,x)\dx
                \\ & \quad 
                =\iint_{\Do} \left(-\frac{1}{2}\wlim{q^2}
                +\left(P-u^2\right)\right)\test \dx\dt,
                \quad \forall \test \in C^\infty_c(\Do).
        \end{split}
\end{equation}
\end{lemma}

\begin{proof}
Starting from \eqref{eq:fq-disc-II} with $f(q)=q$, we argue as 
in the proof of Lemma \ref{lem:weak-f(q)} to conclude 
the validity of \eqref{eq:weak-q}.
\end{proof}

The next lemma tells us that the weak limits in 
Lemma \ref{lem:boundq2} satisfy the 
initial data in an appropriate sense.

\begin{lemma}\label{lem:ex1}
  Let $q$ and $\wlim{q^2}$ be the weak limits identified in Lemma
  \ref{lem:boundq2}. Then
  \begin{equation}
    \label{eq:weaklimits-initialdata}
    \begin{split}
      &\lim_{t\to 0}\int_\R q^2(t,x)\dx =\int_\R (\px u_0)^2\dx,
      \\
      & \lim_{t\to 0}\int_\R \wlim{q^2}(t,x)\dx =\int_\R (\px
      u_0)^2\dx.
    \end{split}
  \end{equation}
\end{lemma}

\begin{proof} 
The proof is similar to that in \cite{Coclite:2006kx}.
\end{proof}

We can now wrap up the proof of the 
strong convergence of $\seq{\qdx}_{\Dx>0}$.

\begin{lemma}
Let $q$ and $\wlim{q^2}$ be the weak limits 
identified in Lemma \ref{lem:boundq2}. Then 
\begin{equation}\label{eq:strongconv-tmp}
        \text{$\wlim{q^2}(t,x)=q^2(t,x)$ for a.e.~$(t,x)\in \Do$.}
\end{equation}
Consequently, as $\Dx\to 0$ (along a subsequence if necessary)
\begin{equation}\label{eq:strongconv}
        \text{$\qdx \to q$ in $L^{2}_{\loc}(\Do)$ and a.e.~in $\Do$}.
\end{equation}
\end{lemma}

\begin{proof}
By Lemma \ref{lem:weak-f(q)},
\begin{equation}\label{eq:eq03-summarize}
        \pt \wlim{f(q)}+\px\left(u \wlim{f(q)}\right)
        \le\wlim{qf(q)}-\frac{1}{2}\wlim{f'(q)q^2}
        +\left(u^2-P\right) \wlim{f'(q)},
\end{equation}
in the sense of distributions on $\Do$, for any convex 
function $f\in C^1(\R)$ with $f'$ bounded
Moreover, by Lemma \ref{lem:weak-q},
\begin{equation}\label{eq:q1-summarize}
        \pt q+\px(u q)=\frac{1}{2}\wlim{q^2}+ u^2-P,
\end{equation}
in the sense of distributions on $\Do$.  Equipped 
with \eqref{eq:eq03-summarize}, \eqref{eq:q1-summarize}, 
\eqref{eq:weaklimits-initialdata}, and \eqref{eq:oleinik_limit}, we can argue 
exactly as in Xin and Zhang \cite{Xin:2000qf} to arrive at 
\eqref{eq:strongconv-tmp}. In view of Lemma \ref {lem:prelim}, claim 
\eqref{eq:strongconv} follows immediately from 
\eqref{eq:strongconv-tmp} and \eqref{eq:udxhigh}.
\end{proof}

We now prove that the limit $u$ satisfies \ref{def:soldistri}.
\begin{lemma}\label{lem:weaksol} 
For any $\test\in C^\infty_c(\Do)$,
\begin{equation}\label{eq:uweaksol}
        \begin{gathered}
                \int_0^T\int_\R u\test_t + \left(\frac{u^2}{2}+P\right)\test_x \dx\dt  = 0,
                \\ \int_0^T\int_\R P\left(\test-\test_{xx}\right)\dx\dt
                =\int_0^T\int_\R \left(u^2
                +\frac{1}{2}\left(\partial_x u\right)^2\right)\test\dx\dt
        \end{gathered}
\end{equation}
\end{lemma}

\begin{proof}
It is not difficult to establish the equation for $P$, since 
we have already established that $\partial_x\udx\to \partial_x u$ 
in $L^2_{\loc}(\Do)$, cf.~\eqref{eq:strongconv}. 
Indeed, we have
\begin{align*}
        \int_0^T& \int_\R \Pdx \left(\test-\test_{xx}\right)\dx\dt \\
        &= \sum_{n,j} \iint_{I^n_{j-1/2}} \Pdx \left(\test-\test_{xx}\right)\dx\dt\\
        &= \sum_{n,j}\Pnj \iint_{I^n_{j-1/2}} \test-\test_{xx}\dx\dt \\
        &\qquad + \underbrace{\sum_{n,j} \iint_{I^n_{j-1/2}}\left(\Pnj-\Pdx\right)
        \left(\test-\test_{xx}\right)\dx\dt}_{w_1}\\
        &= \Dt\Dx \sum_{n,j}\Pnj\left(\test^n_j-\Dm\Dp\test^n_j\right)\\
        &\qquad + w_1 + \underbrace{ \sum_{n,j}\Pnj \iint_{I^n_{j-1/2}}
        \left( \left(\test-\test^n_j\right)
        +\left(\test_{xx}-\Dm\Dp\test^n_j\right)\right) \dx\dt }_{w_2}.
\end{align*}
Since $\abs{\Pdx(t,x)-\Pnj}\le C\Dx$ for $(t,x)\in I^n_{j-1/2}$,
$w_1\to 0$ as $\Dx\to 0$.  Similarly, since $\test-\test_{xx}$ is
close to $\test^n_j-\Dm\Dp\test^n_j$ in $I^n_{j-1/2}$ and $\Pnj$ is
bounded, we conclude that $w_2\to 0$ as $\Dx\to 0$.  Set
\begin{align*}
        f^n_j &= \left(\unjmh\wedge 0\right)^2+ \left(\unjph\vee 0\right)^2
        +\frac{1}{2}\left(\qnj\right)^2,\\
        f_{\Dx}&=\left(\udx\right)^2+\frac{1}{2}\left(\partial_x\udx\right)^2.
\end{align*}
Using the scheme for $\Pnj$, cf.~\eqref{eq:P-disc},
\begin{align*}
  & \Dt\Dx\sum_{n,j} \Pnj\left(\test^n_j-\Dm\Dp\test^n_j\right)
  =\Dt\Dx\sum_{n,j}f^n_j \test^n_j 
  \\ & \qquad =\int_0^T\int_\R f_{\Dx}\test\dx\dt \\
  &\qquad\qquad + \underbrace{\sum_{n,j}
    \iint_{I^n_{j-1/2}}\left(f^n_j-f_{\Dx}\right)\test\dx\dt}_{w_3}
  +\sum_{n,j} f^n_j \underbrace{\iint_{I^n_{j-1/2}}
    \left(\test-\test^n_j\right)\dx\dt}_{=0}.
\end{align*}
By the definition of $\qdx$ we have that 
$\qdx=\qnj + (t-t^n)\DT\qnj$ for $t\in I^n$. Hence
$$
(\qdx)^2-\left(\qnj\right)^2
=2(t-t^n)\qnj\DT\qnj
+\left(t-t^n\right)^2\left(\DT\qnj\right)^2 \quad\text{in $I^n_{j-1/2}$.}
$$
so that, assuming $\mathrm{supp}(\test)\subset [x_{j_a},x_{j_b}]$ 
for some intergers $j_a$ and $j_b$,
\begin{align*}
  &\sum_{n,j} \iint_{I^n_{j-1/2}}
  \abs{\left(\qnj\right)^2-\left(\qdx\right)^2}\abs{\test}\dx\dt \\ 
  & \quad \le C \Dx\sum_n \sum_{j_a}^{j_b}\int_{t^n}^{t^{n+1}}
  \left(2(t-t^n)\qnj\DT\qnj+\left(t-t^n\right)^2\left(\DT\qnj\right)^2\right) \dt
  \\ & \quad = C\Dt\Dx\sum_n\sum_{j=j_a}^{j_b}
  \left(\Dt \qnj\DT\qnj+\left(\Dt\DT\qnj\right)^2\right) \\
  & \quad =\mathcal{O}\left(\sqrt{\Dt\Dx^\theta}+\Dt\Dx^\theta\right).
\end{align*}
By the H\"older continuity of $\udx$; recall that $\udx\in C^{0,\ell}$
with $\ell=1-2/(2+\alpha)$, we find
$$
\abs{\left(\unjmh\wedge 0\right)^2+ \left(\unjph\vee 0\right)^2 -
\left(\udx\right)^2}=\mathcal{O}(\Dx^\ell+\Dt^\ell),
$$
and therefore $w_3\to 0$ as $\Dx\to 0$.  Hence, using \eqref{eq:uuniform},  
\eqref{eq:der1}, and \eqref{eq:strongconv},
\begin{align*}
        \int_0^T\int_\R P\left(\test-\test_{xx}\right)\dx\dt
        & =\lim_{\Dx\to 0} \int_0^T\int_\R\Pdx\left(\test-\test_{xx}\right)\dx\dt
        \\ &=\lim_{\Dx\downarrow 0} \int_0^T\int_\R 
        \left(\left(\udx\right)^2+\frac{1}{2}\left(\partial_x\udx\right)^2\right)\test\dx\dt 
        \\ &=\int_0^T\int_\R \left(u^2+\frac{1}{2}\left(\partial_x u\right)^2\right)\test \dx\dt.
\end{align*}
This means that the second equation in \eqref{eq:uweaksol} holds.

To establish the first equality in \eqref{eq:uweaksol}, we derive a divergence-form 
version of the scheme \eqref{eq:u-disc}.  To this end, introduce the 
functions $f_{\vee}(u)=\frac12 (u\vee 0)^2$ and $f_{\wedge}(u)=\frac12 (u\wedge 0)^2$.  
Observe that $f_{\vee}$ and $f_{\wedge}$ are piecewise $C^2$, and the 
absolute value of the second derivatives are bounded by $1$.  
By the discrete chain rule,
$$
\left(\unjph\vee 0 \right)\Dm\unjph = 
\Dm f_{\vee}(\unjph)+
\mathcal{O}\left(\Dx \left(\Dm \unjph\right)^2\right)
$$
and
$$
\left(\unjph\wedge 0 \right)\Dp\unjph = \Dp f_{\wedge}(\unjph)
+\mathcal{O}\left(\Dx \left(\Dp \unjph\right)^2\right).
$$
Consequently, we can replace \eqref{eq:u-disc} by
\begin{equation}\label{eq:u-disc-cons}
        \begin{split}
                \DT\unjph &+ \Dm f_{\vee}(\unjph) +\Dp f_{\wedge}(\unjph) +\Dp\Pnj
                \\ & =\mathcal{O}\left(\Dx \left\{ \left(\Dm \unjph\right)^2
                +\left(\Dp \unjph\right)^2 \right\}\right).
        \end{split}
\end{equation}
Observe that
\begin{equation}\label{eq:diff-relations}
        D=\frac{\Dm+\Dp}{2}, \qquad \Dx \Dm\Dp = \Dp-\Dm,
        \qquad f_{\vee} +f_{\wedge}= \frac{u^2}{2}.
\end{equation}
Using these identities, we can restate \eqref{eq:u-disc-cons} as
\begin{equation}\label{eq:u-disc-cons-II}
        \begin{split}
                \DT\unjph & + \Dm\left[\frac{\left(\unjph\right)^2}{4}
                +\frac{1}{2}\left(f_{\vee}(\unjph)
                -f_{\wedge}(\unjph)\right)\right]
                \\ &+\Dp \left[\frac{\left(\unjph\right)^2}{4}
                +\frac{1}{2}\left(f_{\wedge}(\unjph)-f_{\vee}(\unjph)\right)\right]+\Dp\Pnj
                \\ &\qquad\qquad = {\mathcal O}\left(\Dx \left\{ \left(\Dm
                \unjph\right)^2 + \left(\Dp \unjph\right)^2\right\}\right).
        \end{split}
\end{equation}
Using, cf.~\eqref{eq:diff-relations},
\begin{align*}
        &\Dm \left( f_{\vee}(\unjph)- f_{\wedge}(\unjph)\right)
        +\Dp\left(f_{\wedge}(\unjph)-f_{\vee}(\unjph)\right) \\ & \qquad
        =\Dx \Dm\Dp f_{\wedge}(\unjph)-\Dx \Dm\Dp f_{\vee}(\unjph),
\end{align*}
equation \eqref{eq:u-disc-cons-II} becomes
\begin{equation}\label{eq:u-disc-cons-III}
        \begin{split}
                \DT\unjph & + D \left(\frac{\left(\unjph\right)^2}{2}\right) +\Dp\Pnj \\
                &=\mathcal{O}\left(\Dx \left\{ \left(\Dm \unjph\right)^2
                +\left(\Dp \unjph\right)^2\right\}\right) \\ & \qquad
                +\Dx\left\{\Dm\Dp f_{\vee}(\unjph)
                -\Dm\Dp f_{\wedge}(\unjph)\right\}.
        \end{split}
\end{equation}
Now fix $\test\in C^2_c(\Do)$ and define $\test^n_j$ as
before, cf.~\eqref{eq:testdef}. 
Multiplying \eqref{eq:u-disc-cons-III} by $\test^n_j\Dt\Dx$ and 
performing partial summations gives
\begin{align*}
        \underbrace{\Dt\Dx\sum_{n,j} \unjph \DT\test^n_j}_{E_1}
        \!+\!\underbrace{\Dt\Dx\sum_{n,j}\frac{\left(\unjph\right)^2}{2}D\test^n_j}_{E_2}
        \!+\!\underbrace{\Dt\Dx\sum_{n,j}\Pnj\Dm\test^n_j}_{E_3}
        =\mathcal{O}(\Dx),
\end{align*}
by using \eqref{eq:enest}. We have $\abs{\udx-\unjmh}
\le C\left(\Dx\abs{\qnj}+\Dt\abs{\DT\unjph}\right)$. 
Using this and \eqref{eq:uuniform}, we compute as follows:
\begin{align*}
        E_1&=\int_0^T\int_\R \udx\test_t\dx\dt
        +\sum_{n,j}\iint_{I^n_{j}} \left(\unjph-\udx\right)\test_t\dx\dt\\
        &\qquad\qquad+\sum_{n,j}\unjph\iint_{I^n_{j}}\left(\DT\test^n_j-\test_t\right)\dx\dt\\
        &=\int_0^T\int_\R \udx\test_t\dx\dt + \mathcal{O}(\Dx)
        \to \int_0^T\int_\R u\test_t\dx\dt\quad\text{as $\Dx\to 0$.}
\end{align*}
In the same way, equipped with \eqref{eq:uuniform} and 
\eqref{eq:Pconvergence}, we can show that
$$
E_2 \to \int_0^T\int_\R \frac{u^2}{2}\test_x\dx\dt, 
\quad E_3 \to \int_0^T\int_\R P\test_x\dx\dt
\quad \text{as $\Dx\to 0$.}
$$
thus proving the first equality in \eqref{eq:uweaksol}. 
This concludes the proof of the lemma.
\end{proof}

\section{Numerical examples}\label{sec:numex}
We have tried the difference method presented here on several
examples, and in doing this found that the first order method analyzed
in this paper exhibits very slow convergence, and thus requires a very
small mesh size $\Dx$ to compute reasonable solutions.  This is not
surprising and appears to be the case with other schemes in the
literature as well.  Therefore we have implemented a second order
extension of the method. This second order extension is based on the
conservative version of the scheme
\begin{multline}
  \DT \unjph + \Dm\left[
    \left(\unjph\vee 0\right) \unjph + \left(\unjph\wedge
      0\right)u^n_{j+3/2} + P^n_{j+1}\right] \\
  =\Dm\left(\unjph\vee 0\right)\unjmh + \Dm\left(\unjph\wedge 0\right)\unjph,
\end{multline}
which can be viewed as a balance equation with a flux 
across $x=x_{j+1/2}$ given by 
$$
F_{j+1/2}^n=\left(\unjph\vee 0\right) \unjph + \left(\unjph\wedge
      0\right)u^n_{j+3/2} + P^n_{j+1}.
$$
Taking this viewpoint, we \emph{define} the second order finite volume
scheme by 
\begin{equation}
  \label{eq:second-disc}
  \DT u_j^n + \Dm F_{j+1/2}^{n+1/2} = u_j^n 
  \Dm\left(u^{n+1/2}_{j+1/2}\right). 
\end{equation}
Here $u^{n+1/2}_{j+1/2}$ is a first order approximation of the value
at the point $x=x_{j+1/2}$, $t=t^n+\Dt/2$. 
This
approximation is found by setting 
$$
\unjph = \frac{1}{2}\left(u^n_{j}+u^n_{j+1}\right),
$$
and then using the scheme \eqref{eq:u-disc} for half a time step
(i.e., $\Dt/2$). This scheme is  a formally second order accurate
finite volume approximation, and this simple adaptation produces
significantly more accurate approximations. 

\begin{figure}[htbp]
  \centering
  \includegraphics[width=0.6\linewidth]{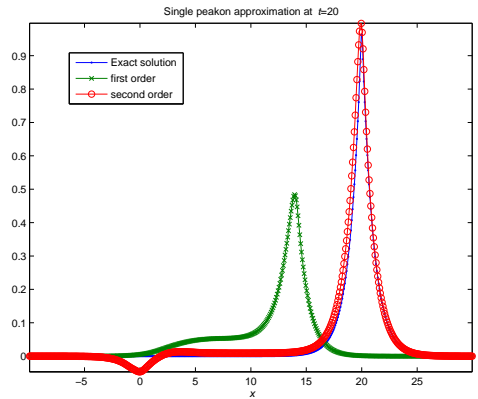}
  \caption{Approximations using $\Dx=40/2^9$, at $t=20$ to the single peakon.}
  \label{fig:1}
\end{figure}

In Figure~\ref{fig:1} we show the approximations calculated by the
first order scheme \eqref{eq:u-disc} and the second order scheme
\eqref{eq:second-disc} for the single peakon example. In this case the
exact solution reads
$$
u(x,t)=e^{-\abs{x-t}}.
$$
Figure~\ref{fig:1} shows the solutions calculated using $2^9$ equally
spaced grid points in the interval $[-10,30]$ for $t=20$. We see that
the second order method is much more accurate than the first order
method. 
In passing, we note that we have not used the strict CFL-condition
\eqref{eq:cfl}, but the more natural condition
$$
\Dt \le \max_j\seq{u^n_j} \Dx. 
$$
This holds for the second order scheme as well. 

In order to investigate the convergence properties of the two methods,
we computed errors in $L^1$ for the two schemes. Table~\ref{tab:1}
shows the computed $L^1$ errors in the case of a single peakon under
mesh refinement. In this context, the $L^1$ error is defined as
$$
\text{$L^1$ error} = \Dx\sum_i \abs{u_{\Dx}(x_i,t) - u(x_i,t)},
$$
where $u$ is the exact solution. We used $t=20$ 
and $\Dx=40/2^k$ for $k=5,\ldots,12$. As expected,
and as reported in \cite{Coclite:2006kx},
the first order method converges very slowly.
\begin{table}[htbp]
  \centering
  \begin{tabular}[h]{l|rrrrrrrrr}
    $k$ & 5 & 6 & 7 & 8 & 9 & 10 & 11 & 12 & 13 \\
    \hline
    $1^\mathrm{st}$ &   2.92    &3.23   &3.41   &3.53   &3.57   &3.51   &3.32   &3.01   &2.64\\
    $2^\mathrm{nd}$ & 5.36      &5.17   &3.29   &1.27   &0.60   &0.36   &0.21   &0.13   &0.09 \\
  \end{tabular}
  \caption{$L^1$ errors for the single peakon case, at $t=20$, for
    $x\in [-10,30]$, $\Dx=40/2^k$, $k=5,\ldots,13$}
  \label{tab:1}
\end{table}
One other notable feature of Table~\ref{tab:1} is that the second
order method seems to converge at a rate slightly less than $1$. 

\begin{figure}[htbp]
  \centering
  \begin{tabular}{lr}
    \includegraphics[width=0.5\linewidth]{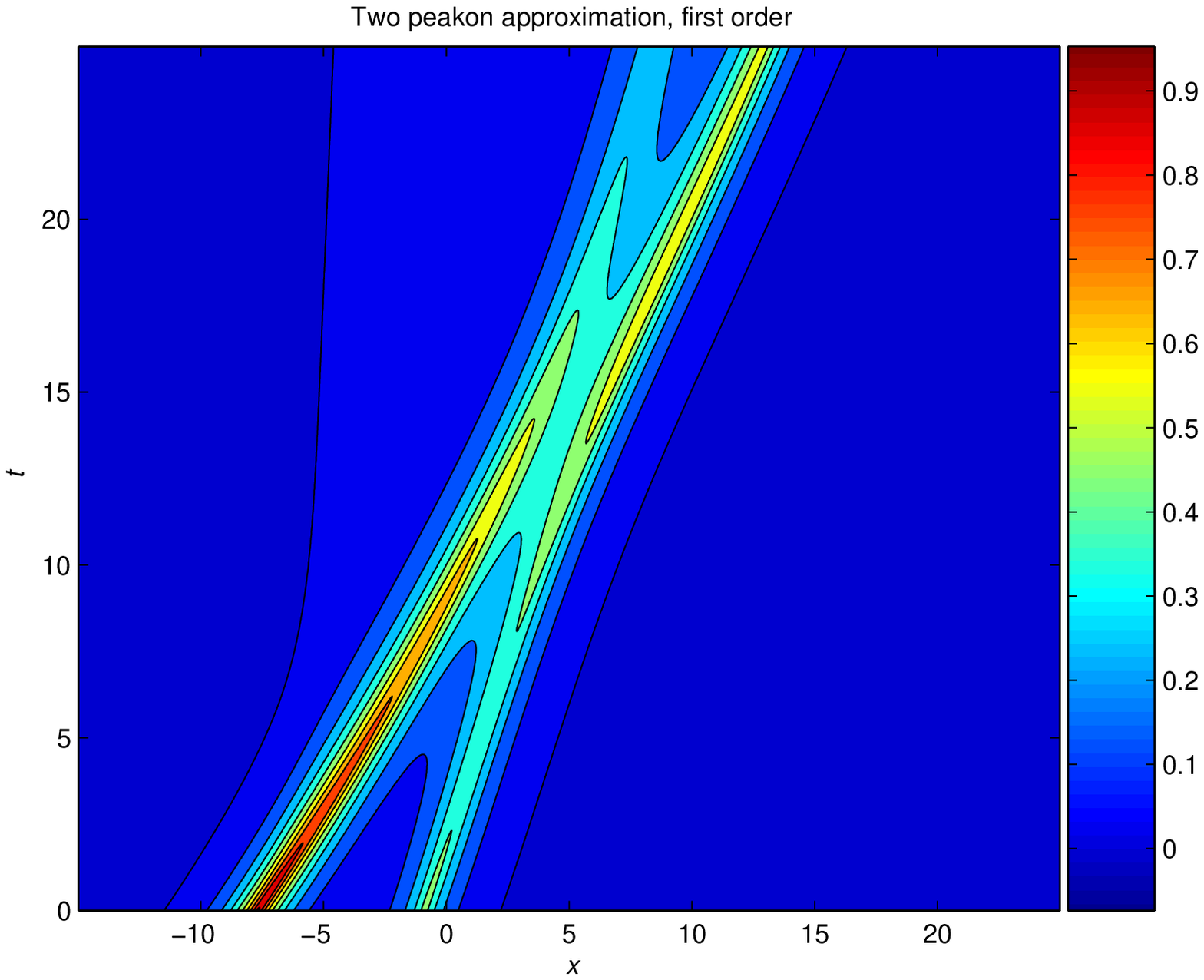}&
    \includegraphics[width=0.5\linewidth]{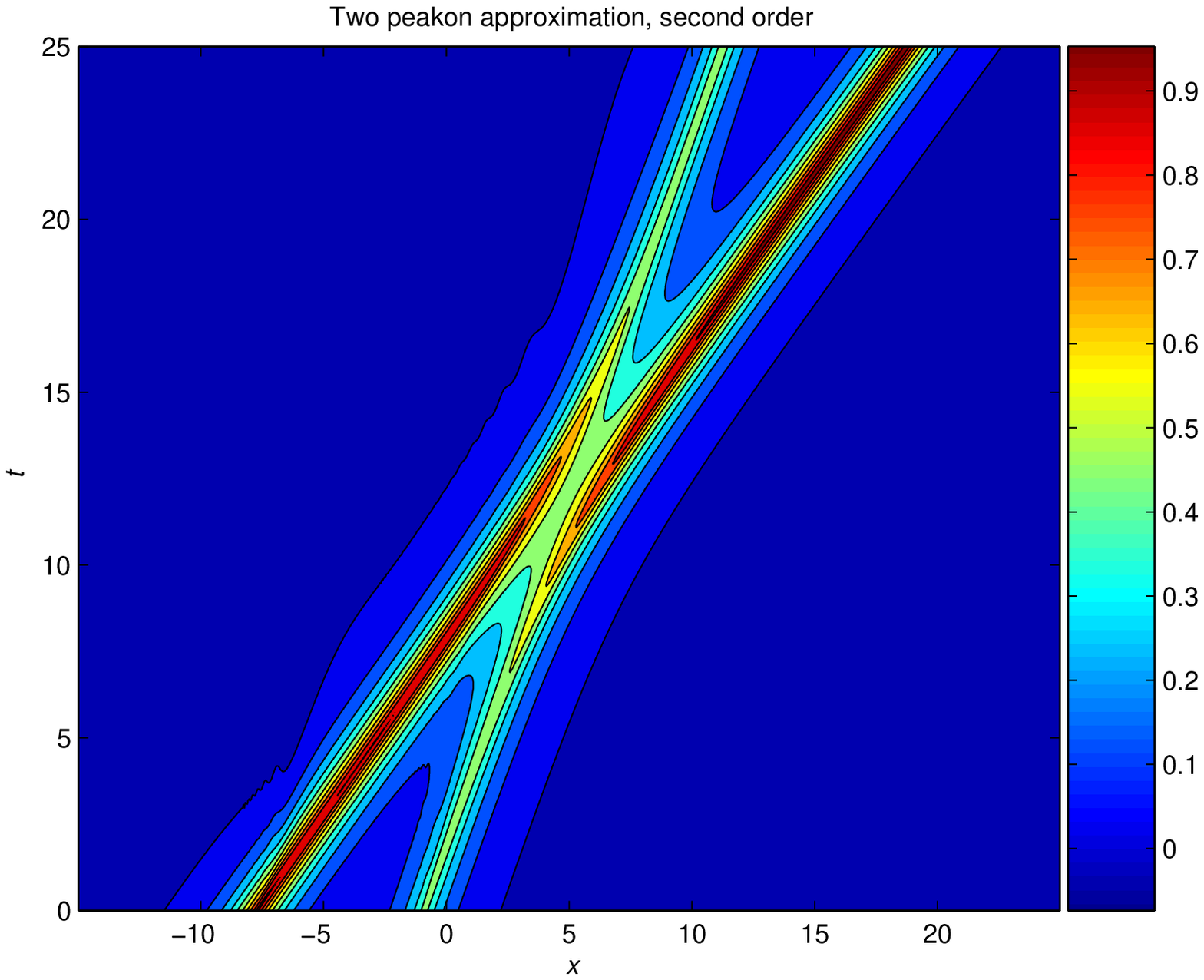}
  \end{tabular}
  \caption{Approximations to \eqref{eq:twopeak} using
    $\Dx=40/2^{10}$. Left: first order method\eqref{eq:u-disc}. Right :
    second order method \eqref{eq:second-disc}.} 
  \label{fig:2}
\end{figure}

\begin{figure}[htbp]
  \centering
  \includegraphics[width=0.6\linewidth]{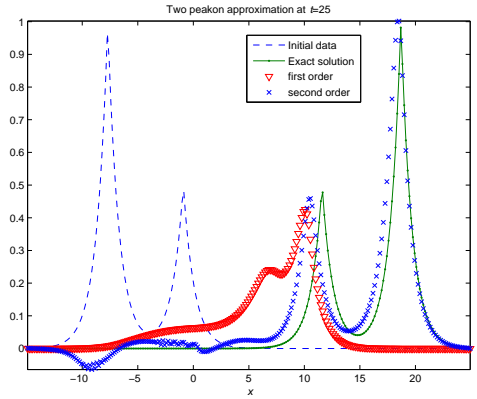}
  \caption{The approximations to \eqref{eq:twopeak} at $t=25$ and $\Dx=40/2^8$.}
  \label{fig:3}
\end{figure}

The two-peakon solution is considerably more complicated than the
single peakon, and this is also a much harder challenge
computationally, see e.g.,~\cite{Artebrant:2006qy} and
\cite{Kalisch:2006kx}. We use the two-peakon solution given by
\begin{equation}
  \label{eq:twopeak}
  u(x,t)=m_1(t)e^{-\abs{x-x_1(t)}} + m_2(t)e^{-\abs{x-x_2(t)}},
\end{equation}
with
$$
\begin{gathered}
  x_1(t)=\log\left(\frac{18e^{t-10}}{e^{(t-10)/2} + 6}\right),\quad 
  x_2(t)=\log\left(40e^{t-10}+60e^{(t-10)/2}\right)\\
  m_1(t)=\frac{e^{(t-10)/2}+6}{2e^{(t-10)/2}+3}, \quad
  m_2(t)=\frac{e^{(t-10)/2}+\frac{2}{3}}{e^{(t-10)/2}+3}.
\end{gathered}
$$
These formulas were taken from
\cite{Lundmark:2005kz}. Figure~\ref{fig:2} shows a contour plot of the
approximate solutions found by using the first and second order
methods, and $\Dx=40/2^{10}$ for $x\in [-15,25]$ and $t\in [0,25]$. 
We see that the interaction between the two peakons is poorly
represented by the first order method. Both the location as well as
the magnitude of the peaks are far from the correct value. This is
also illustrated Figure~\ref{fig:3} where we show the approximations using
$\Dx=40/2^{8}$ at $t=25$.
We have also calculated errors for the two-peakon case. Indeed for
$\Dx\ge 40/2^{12}$, the first order method did not seem to converge,
and in order to give meaningful answers, this method demands very fine
discretizations. These results are reported in Table~\ref{tab:2}. For
$\Dx > 40/2^8$ none of the methods gave satisfactory results.
\begin{table}[htbp]
  \centering
  \begin{tabular}[h]{l|rrrrrr}
    $k$&  8&  9& 10& 11& 12 & 13\\
    \hline
    $1^{\text{st}}$ &4.56&3.64&3.97&4.18&4.05 & 3.70\\
    $2^{\text{nd}}$ &1.88&1.04&0.63&0.38&0.22 & 0.16
  \end{tabular}
  \caption{$L^1$ errors for the approximation to
    \eqref{eq:second-disc}, $t=25$, $x\in [-15,25]$, $\Dx=40/2^k$, $k=8,\ldots,13$. } 
  \label{tab:2}
\end{table}

\begin{figure}[htbp]
   \centering
   \begin{tabular}{lr}
     \includegraphics[width=0.5\linewidth]{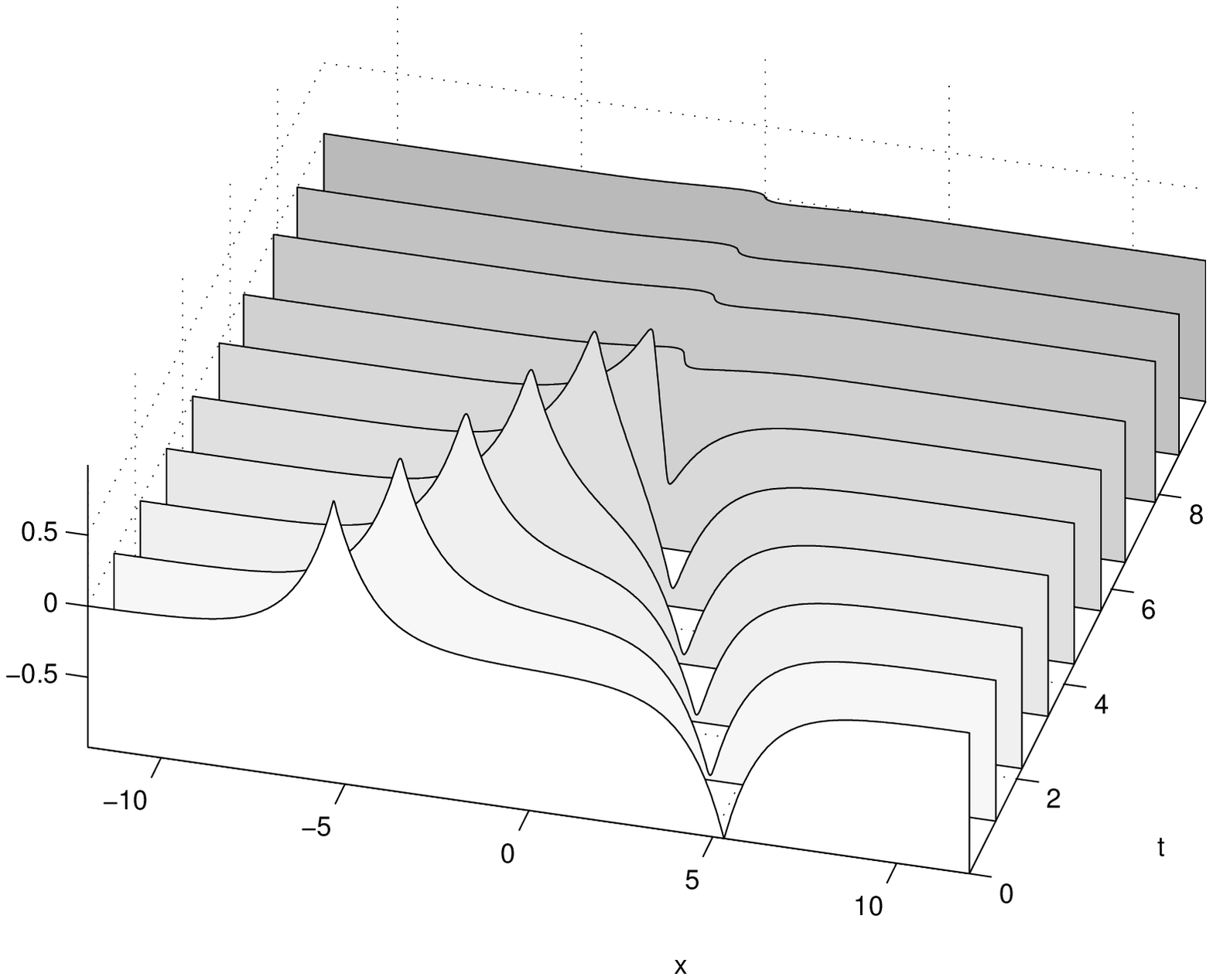}&
     \includegraphics[width=0.5\linewidth]{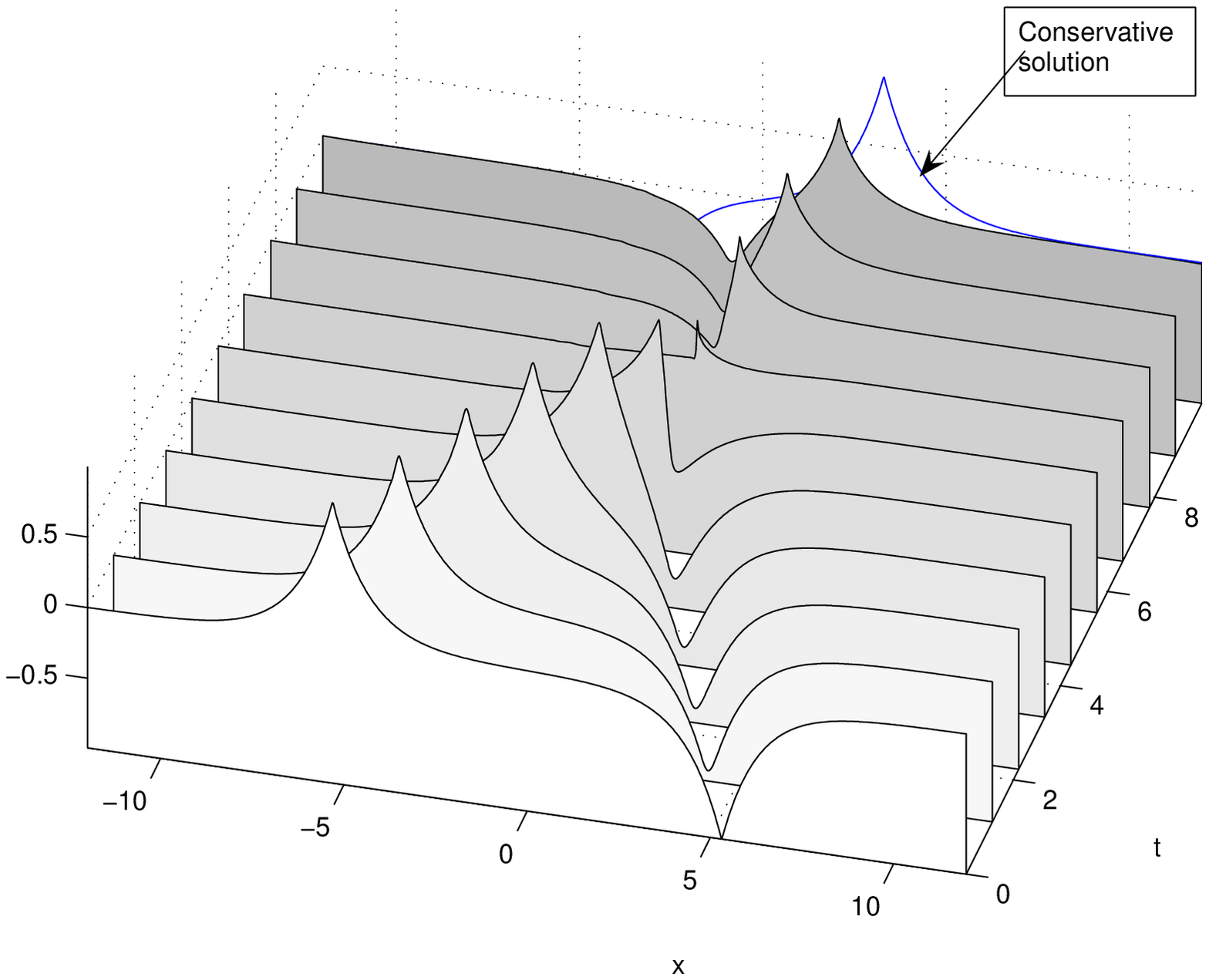}
   \end{tabular}
   \caption{The numerical solutions to the initial value problem 
     \eqref{eq:peakonantipeakon}. Left: first order method, right:
     second order version.}
   \label{fig:peakantipeak}
 \end{figure}

In our final example we choose initial data corresponding 
to a peakon-antipeakon collision:
\begin{equation}
  \label{eq:peakonantipeakon}
  u_0(x)=-\tanh(6)\left(e^{-\abs{x+y(6)}}-e^{-\abs{x-y(6)}}\right),
\end{equation}
where $y(t)=\log(\cosh(t))$. 
In this case we have a ``peakon anti-peakon collision'' at $t=6$.  In
Figure~\ref{fig:peakantipeak} we exhibit the approximations generated
by the first order (left) and the second order method for $t\in [0,10]$
and $\Dx=24/2^{12}$.  It is
clear that the first order scheme generates the dissipative solution, and for $t$
larger than the collision time, the first order approximation
vanishes. Regarding the second order approximation, it seems to
continue as a peakon moving to the right, and an anti-peakon moving to
the left. The magnitudes and speeds of these features are however far
from the conservative solution, and we have indicated the conservative
solution in the right hand figure. 


\end{document}